\documentclass[11pt,reqno]{amsart}
\usepackage{graphicx,tikz}
\usepackage{amsfonts,amsmath,amssymb}
\usetikzlibrary{matrix,arrows}
\newcommand{\rl}{{\mathbb{R}}}
\newcommand{\cx}{{\mathbb{C}}}

\newcommand{\G}{\Gamma}

\newcommand{\e}{\varepsilon}

\newcommand{\tmop}[1]{\ensuremath{\operatorname{#1}}}
\renewcommand{\Re}{\tmop{Re}}
\renewcommand{\Im}{\tmop{Im}}

\newcommand{\B}{\mathbb{B}}
\newcommand{\D}{\mathbb D}
\newcommand{\C}{\mathbb C}
\newcommand{\HH}{\mathbb H}
\newcommand{\R}{\mathbb R}

\newcommand{\PP}{\mathbb P}

\newtheorem{theorem}{Theorem}[section]
\newtheorem{lemma}[theorem]{Lemma}
\newtheorem{prop}[theorem]{Proposition}
\newtheorem{cor}[theorem]{Corollary}

\theoremstyle{definition}

\begin{document}
\title{Some aspects of holomorphic mappings: a survey }
\author{Sergey Pinchuk, Rasul Shafikov and Alexandre Sukhov}
\begin{abstract}
This expository paper is concerned with the properties of proper holomorphic mappings between domains in complex 
affine spaces. We discuss some of the main geometric methods of this theory, such as  the Reflection Principle, 
the scaling method, and the Kobayashi-Royden metric. We sketch the proofs of certain principal results and discuss some 
recent achievements. Several open problems are also stated.
\end{abstract}

\maketitle

\let\thefootnote\relax

* Department of Mathematics, Indiana University, 831 E 3rd St. Rawles Hall, Bloomington,   IN 47405, USA, e-mail: pinchuk@indiana.edu

** Department of Mathematics, the University of Western Ontario, London, Ontario, N6A 5B7, Canada,
e-mail: shafikov@uwo.ca. The author is partially supported by the Natural Sciences and Engineering 
Research Council of Canada.

*** Universit\'e  de Lille (Sciences et Technologies), 
U.F.R. de Math\'ematiques, 59655 Villeneuve d'Ascq, Cedex, France,
e-mail: sukhov@math.univ-lille1.fr. The author is partially supported by Labex CEMPI.

\tableofcontents

\section{Introduction}
This expository paper is dedicated to geometric properties of holomorphic mappings 
between domains in complex affine spaces (in general of different dimensions). The first results in this direction (mainly in 
complex dimension 2) are due to H.~Poincar\'e, E.~Cartan and B.~Segre. The rigidity of complex structures with boundary--one of the main phenomena of complex analysis in higher dimensions--was already discovered and studied in these classical works. The next major step in this theory was made in the 70-ties with intensive investigation of the geometry of strictly pseudoconvex domains. Further progress concerns more general classes of domains (weakly pseudoconvex or not pseudoconvex at all). In this survey paper we try to present  some of the main ideas in the development of the theory. Our presentation is certainly incomplete, as, regrettably, many important topics and results were not included in the scope of the paper. An interested reader may become acquainted  with them using other monographs  and expository papers \cite{BaEbRo1, BaEbRo2, Bed2, Ber6, Ch3, DiPi5, For6, For6b, For8, GKK, Is, IsKr, Kr, Me2, Me3, Pi10, Si3, Tu2, Vi, Vi2}.

The authors are deeply grateful to referees for many corrections and suggestions.

\section{Preliminaries}

 Denote by $z = (z_1,...,z_n)$ the standard complex coordinates in $\C^n$. We often use the (vector) notation $z = x +iy$ for the real and imaginary parts. Denote by $\vert z \vert$ the euclidean norm of $z$ and by $(z,w) = \sum_j z_j \overline{w_j}$ the hermitian inner product. We also use the notation $\langle z,  w \rangle = (z,\overline w)$. 
 
 As usual, a {\it domain} $\Omega$ in $\C^n$ is a connected open subset of $\C^n$. Denote by $\partial\Omega$ the boundary 
 of~$\Omega$.
The unit ball of $\C^n$ is denoted by  $\B^n = \{ z \in \C^n : \vert z \vert < 1 \}$, while for $n = 1$ we use the notation $\D := \B^1$ for the unit disc in $\C$. The ball $p+ r\B^n$ of radius $r > 0$ centred at a point $p \in \C^n$ will be denoted by $\B^n(p,r)$. Another basic example of a domain in $\mathbb C^n$ is the unit polydisc $\D^n$, or more generally, $\D^n(p,r):= p + r\D^n$. 
Finally,  
\begin{equation}\label{e.H}
\HH = \{ z \in \C^n: 2\Re z_n + \vert z_1 \vert^2 +...+ \vert z_{n-1} \vert^2 < 0 \}
\end{equation} 
is an unbounded realization of the unit ball $\B^n$.

\subsection{Classes of functions}  Denote by ${\mathcal O}(\Omega)$ the class of holomorphic functions in a domain $\Omega$. If $\Omega'$ is a domain in $\C^m$, we use the notation 
${\mathcal O}(\Omega,\Omega')$ for the class of holomorphic mappings from 
$\Omega$ to $\Omega'$.

For a positive integer $k$, $C^k(\Omega)$ denotes the space of $C^k$-smooth complex-valued functions in~$\Omega$. Also $C^k(\overline\Omega)$ denotes the class of functions whose partial derivatives up to order $k$ extend as continuous functions on 
$\overline\Omega$. If $s > 0$ is a real noninteger and $k$ is its integer part, $C^s(\Omega)$ denotes the space of functions of class $C^k(\Omega)$ such that their partial derivatives of order $k$ are (global) H\"older-continuous in $\Omega$ with the exponent $s-k$; these derivatives  automatically satisfy the  H\"older condition on $\overline\Omega$ so the notation $C^s(\overline\Omega)$ for the same space of functions is also appropriate. Finally, we denote by $PSH(\Omega)$ the class of plurisubharmonic functions on $\Omega$.

Let $\Omega$ be a domain in $\C^n$ and $f:\Omega \to \C^N$ be a vector function (not necessarily holomorphic or smooth) on $\Omega$. Let $\gamma$ be a subset of the boundary $\partial\Omega$. The {\it cluster set} $C_\Omega(f;\gamma)$ of $f$ on $\gamma$ is defined as the set of all limit points of the sequences $\{f(z^k)\}$ in $\C^N$, where $\{z^k\}$ is any sequence 
in~$\Omega$ converging to a point in~$\gamma$. The cluster set $C_\Omega(f;\gamma)$ is empty if $\lim \vert f(z) \vert = +\infty$ when $z \to\gamma$. Note that a holomorphic map $f:\Omega \to \Omega'$ between two domains is proper if and only if the cluster set $C_\Omega(f;\partial\Omega)$ does not intersect $\Omega'$. For bounded domains one can state this property in the equivalent form: $C_\Omega(f;\partial\Omega) \subset \partial\Omega'$.

\subsection{Real submanifolds of complex spaces.} A (closed) real submanifold $E$ of a domain $\Omega \subset \C^n$ is of class $C^s$ (resp. real analytic) if for every point $p \in E$ there exists an open neighbourhood $U$ of $p$ and a map $\rho: U \longrightarrow \R^d$ of the maximal rank $d<2n$ and of class $C^s$ (resp. real analytic) such that $E \cap U = \rho^{-1}(0)$; then $\rho$ is called a local defining (vector) function of $E$. The positive integer $d$ is the real codimension of $E$. In the fundamental special case $d=1$ we obtain the class of real hypersurfaces.

Let $J$ denote the standard complex structure of $\C^n$. In other words, $J$ acts on a vector $V$ by multiplication by $i$. For every $p \in E$ the {\it holomorphic tangent space}  $H_pE:= T_pE \cap J(T_pE)$ is the maximal complex subspace of the tangent space $T_pE$ of $E$ at $p$. Clearly $H_pE =   \{ V \in \C^n: \partial \rho(p) V = 0 \}$. The complex dimension of $H_pE$ is called the CR dimension of $E$ at $p$; a manifold $E$ is called a {\it CR  (Cauchy-Riemann) manifold} if its CR dimension is independent of $p \in E$.

A real submanifold $E \subset \Omega$ is called {\it generic} (or generating) if the complex span of $T_pM$ coincides with $\C^n$ for all $p \in E$. Note that every generic manifold of real codimension $d$ is a CR manifold of CR dimension $n-d$. A function $\rho = (\rho_1,...,\rho_d)$ defines a generic manifold if $\partial\rho_1 \wedge ...\wedge \rho_d \neq 0$.
Of special importance are the so-called {\it totally real manifolds}, i.e., submanifolds $E$ for which $H_pE = \{ 0 \}$ at every $p \in E$. A totally real manifold is generic if and only if its real dimension is equal to $n$; this is the maximal possible value for the dimension of a totally real manifold. 

Let $E$ be a generic manifold of real codimension $d$ contained in the boundary $\partial\Omega$ of a domain $\Omega$ in $\C^n$. Our considerations are local. Consider tangent vector fields $X_j$, $j = 1,...,n-d$, on $E$ (of type (1,0))  which form a basis in the space of local sections of the holomorphic tangent bundle $H(E)$ near $p$. A $C^1$-smooth function $f$ on $E$ is called a CR function if it satisfies the first order PDE system on $E$
\begin{eqnarray}
\label{tangentCR}
X_jf = 0, \,\, j=1,...,n-d. 
\end{eqnarray}  
These are the {\it tangential Cauchy-Riemann equations}. By  Stokes' formula the equations (\ref{tangentCR}) can be rewritten in 
the equivalent form 
 $[E](f\overline\partial\phi) = 0$ for every test $(n,n-d)$ form $\phi$ on $E$; here $[E]$ denotes the current of integration over $E$. In this weak formulation the notion of a CR function can be extended to the class of continuous or locally integrable functions on $E$. 
 If $E$ is a hypersurface ($d=1$) given by a defining function $\rho$ with $\partial \rho/\partial z_n \neq 0$, then we may choose
$$
X_j = \frac{\partial \rho}{\partial z_n}\frac{\partial }{\partial z_j} -\frac{\partial \rho}{\partial z_j}\frac{\partial }{\partial z_n}, \,\, j=1,...,n-1 .$$

The following approximation theorem is due to Baouendi - Tr\`eves \cite{BaTr}:
\begin{theorem}
\label{BaTr}
Let $M$ be a smooth generic manifold in $\C^n$ and $E \subset M$ be a smooth totally real manifold of dimension $n$. Then 
in a neighbourhood of any point $p \in E$, any CR function $f$ of class $C^s$, $s \ge 0$, on $M$ can be approximated in the 
$C^s$ norm on $M$ by the sequence of holomorphic functions 
$$ ({\bf 1}_E f) * \exp(-k\langle z, z \rangle), \ k = 1,2,...,$$
where ${\bf 1}_E$ denotes the characteristic function of $E$ and the asterisk $*$ denotes the convolution operator.
\end{theorem}
 
\subsection{Pseudoconvex and strictly pseudoconvex domains}   Let $\Omega$ be a bounded domain in~$\C^n$. Suppose that its boundary $\partial\Omega$ is a (compact) real hypersurface of class $C^s$ in $\C^n$. Then there exists a $C^s$-smooth real 
function $\rho$ in a neighbourhood $U$ of the closure $\overline\Omega$ such that $\Omega = \{ \rho < 0 \}$ and 
$d\rho|_{\partial\Omega} \ne 0$. We call such a function $\rho$ a global defining function. If $s \ge 2$ one may consider 
{\it the Levi form} of $\rho$:
\begin{eqnarray}
\label{Leviform}L(\rho,p,V) = \sum_{j,k = 1}^n \frac{\partial^2\rho}{\partial z_j\partial\overline{z}_k}(p)V_j \overline V_k.
\end{eqnarray}
A bounded domain $\Omega$ with $C^2$ boundary is called {\it pseudoconvex} (resp. {\it strictly pseudoconvex}) 
if $L(\rho,p,V) \ge 0$ (resp. $> 0$) for every $V \in H_p(\partial\Omega)$ (resp. every nonzero $V \in H_p(\partial\Omega)$). 
This definition is equivalent to the general notion of pseudoconvexity in the sense of Grauert-Oka: $\Omega$ is pseudoconvex 
if and only if it can be exhausted by a sequence of strictly pseudoconvex domains. Every strictly pseudoconvex domain $\Omega$ admits a global defining function which is 
strictly plurisubharmonic on a neighbourhood $U$ of $\overline\Omega$. The analog of this property for pseudoconvex domains was established by 
Diederich-Fornaess \cite{DiFo1}; 

\begin{theorem}\label{t.df}
Let $\Omega$ be a bounded pseudoconvex domain with $C^s$-boundary, $s \ge 2$. Then there exist 
a $C^s$-smooth defining function $\rho$ in a neighbourhood $U$ of $\overline\Omega$  and a positive $\eta_0<1$ such that for any 
$0 < \eta < \eta_0$, the function $\hat\rho:=-(-\rho)^\eta$ is a strictly 
plurisubharmonic bounded exhaustion function for $\Omega$ 
(i.e., $\hat\rho:\Omega \to (0,a)$ is a proper map for some $a > 0$).
\end{theorem}

The famous example of the so-called ``worm" domain due to the same authors~\cite{DF0} shows that there 
exist smoothly bounded pseudoconvex domains without a plurisubharmonic defining function.

Let $\Gamma$ be a real hypersurface of class $C^2$ in $\C^n$. One can view every holomorphic tangent space $H_p\Gamma$ as an element of the Grassmanian $G(n-1,n)$ of hyperplanes in $\C^n$. Then the holomorphic tangent bundle $H(\Gamma)$  can be viewed as a real submanifold of dimension $2n-1$  of the complex manifold $\C^n \times G(n-1,n)$ of complex dimension $2n-1$. We call it the {\it projectivization of the holomorphic tangent bundle} and denote by $\PP H(\Gamma)$. The following statement, due to Webster \cite{We3}, is easy to check in local coordinates.
\begin{lemma}
$\Gamma$ has a nondegenerate Levi form if and only if $\PP H(\Gamma)$ is a totally real manifold in $\C^n \times G(n-1,n)$.
\end{lemma}

\subsection{Kobayashi-Royden pseudometric} Let $z$ be a point of a domain $\Omega$ and $V$ be a tangent vector at $z$. The infinitesimal Kobayashi-Royden pseudometric $F_\Omega(z,V)$ (the ``length'' of the vector $V$) is defined as 
\begin{equation}
F_\Omega(z,V) = \inf \left\{\lambda>0: \exists\, h \in \mathcal O (\mathbb D, \Omega) {\rm \ with \ } 
h(0)=z,\ h'(0)=\frac{V}{\alpha} \right\} .
\end{equation}
This is a nonnegative upper semicontinuous function on the tangent bundle of $\Omega$; its integrated form coincides with the usual Kobayashi distance. The Kobayashi-Royden metric is decreasing under holomorphic mappings: if $f: \Omega \to \Omega'$ is a holomorphic mapping between two domains in $\C^n$ and $\C^m$ respectively, then
\begin{eqnarray}
\label{Kob1}
F_{\Omega'}(f(z),df(z)V) \le F_\Omega(z,V) .
\end{eqnarray}
In fact, this is the largest metric in the class of (properly normalized) infinitesimal metrics that are decreasing under holomorphic mappings.
It is easy to obtain an upper bound on $F_\Omega$. Indeed, let $\B^n(z,R)$ with $R = {\rm dist\,}(z,\partial\Omega)$ be the ball
contained in $\Omega$. It follows by the holomorphic decreasing property applied to the natural inclusion 
$\iota: \B^n(z,R) \to \Omega$ that the Kobayashi-Royden metric of this ball is bigger than $F_\Omega$. This gives the upper bound
\begin{eqnarray}
\label{Kob2}
F_\Omega(z,V) \le \frac{C \vert V \vert}{{\rm dist\,} (z,\partial\Omega)} .
\end{eqnarray}
Lower bounds require considerably more subtle analysis. Some general estimates can be obtained using plurisubharmonic functions. Sibony \cite{Si1} proposed the  approach based on the following Schwarz-type lemma for subharmonic functions.  For a domain 
$\Omega$ and $z\in \Omega$, denote by $S_z(\Omega)$ the class of functions $u: \Omega \to [0,1]$ such that $u(z) = 0$, $u$ is of class $C^2$ in a neighbourhood of $z$ and $\log u$ is a plurisubharmonic function in $\Omega$.

\begin{lemma}
\label{SiSch}
Let $u\in S_0(\D)$. Then
\begin{itemize}
\item[(a)] $u(\zeta) \le \vert \zeta \vert^2$ for $\zeta \in \D$. The equality holds at some point different from $0$ iff $u(\zeta)$ is identically equal to $\vert \zeta \vert^2$.
\item[(b)] $\Delta u(0) \le 4$ with equality iff $u(\zeta) = \vert \zeta \vert^2$ for every $\zeta \in \D$.
\end{itemize}
\end{lemma}

Consider an infinitesimal pseudometric $P_\Omega$ defined by
\begin{eqnarray}
\label{Kob3}
P_\Omega(z,V) = \sup \{ L(u,z,V)^{1/2}: u \in S_z(\Omega) \} .
\end{eqnarray}
This pseudometric is locally bounded on the tangent bundle by Lemma \ref{SiSch}(a), and is decreasing under holomorphic mappings; hence $$P_\Omega \le F_\Omega.$$
To obtain the estimate from below for Sibony's metric it suffices to construct a function $u \in S_z(\Omega)$ with controlled Levi form. For example, this leads to the following 

\begin{prop}
\label{KobProp1}
Let $\Omega$ be a domain in $\C^n$ and $\rho$ be a negative $C^2$-smooth plurisubharmonic function in $\Omega$. Suppose that the partial derivatives of $\rho$ are bounded on $\Omega$ and there exists a constant $C_1 > 0$ such that 
\begin{eqnarray}
\label{Kob4_1}
L(\rho,z,V) \ge C_1 \vert V \vert^2
\end{eqnarray}
for all $z$ and $V$. Then there exists a constant $C_2 > 0$, depending only on the $C^2$-norm of $\rho$, such that 
\begin{eqnarray}
\label{Kob4}
P_\Omega(z,V) \ge C_2 \left ( C_1^2\frac{\vert \langle \partial \rho(z), V \rangle \vert^2}{\vert \rho(z) \vert^2} + C_1\frac{\vert V \vert^2}{\vert \rho(z) \vert^2} \right ).
\end{eqnarray}
\end{prop}
Note that $\rho$ is not assumed to be a defining function of $\Omega$, although this special case is particularly important in applications. The original argument of Sibony assumes that $\Omega$ is globally bounded but this condition can be dropped. In fact, the estimate (\ref{Kob4}) holds on an open subset of $\Omega$ where (\ref{Kob4_1}) is satisfied. Therefore, it can be used in order to localize the Kobayashi-Royden metric.  Note also that $\Omega$ is not assumed to be bounded or hyperbolic (see \cite{Ber1,Ber2,Su1,Su2, ChCoSu}). In particular, this leads to the following result (see \cite{ChCoSu}):

\begin{prop}
\label{CCS}
Let $\Omega$ be a domain in $\cx^n$, $\rho$ be a plurisubharmonic function in $\Omega$ with $E = \rho^{-1}(0)$, and let $f:\D \longrightarrow \Omega^+ = \{ \rho \ge 0 \}$ be a bounded holomorphic mapping such that the cluster set $C_\D(f,\gamma)$ on an open arc $\gamma \subset \partial \D$ is contained in $E$. Assume that for a certain point $\zeta \in \gamma$ the cluster set 
$C_\D(f,\zeta)$ contains a point $p \in E$ such that, for some $\e > 0$, the function $\rho(z) - \e\vert z \vert^2$ is plurisubharmonic in a neighbourhood of $p$. Then $f$ extends to a H\"older $1/2$-continuous mapping in a neighbourhood of $\zeta$ in $\D \cup \gamma$.
\end{prop}

The proof is based on the estimate (\ref{Kob4}) in a tube neighbourhood of $E$ of the form $\rho < \delta$ with small $\delta > 0$. A special case useful for applications arises when $E$ is a totally real manifold: indeed, such a manifold can be represented as the zero locus of a nonnegative strictly plurisubharmonic function, see \cite{HW}.

\subsection{Some properties of holomorphic functions near real manifolds} Analytic discs form an important special class of holomorphic mappings. Recall that an {\it analytic (or holomorphic) disc} in $\C^n$ is a holomorphic mapping $f:\D \to \C^n$. 
The most interesting case arises when   analytic discs have some boundary regularity (at least, are continuous on $\overline\D$). 
The restriction $f:\partial\D \to \C^n$ is called the boundary of the analytic disc $f$. We say that a disc $f$ is attached or glued 
to a subset $K$ of $\C^n$ if $f(\partial\D) \subset K$.

Let $E$ be a generic submanifold in a domain $\Omega \subset \C^n$ defined as $\{\rho = (\rho_1, ...,
\rho_d) = 0\}$. The {\it wedge} $W(\Omega,E)$ in $\Omega$ {\it with the edge} $E$ is the domain 
$$W(\Omega,E) = \{ z \in \Omega: \rho_j(z) < 0, j=1,...,d \} .$$ 
One can also consider a more general class of domains if we fix an open (convex) cone $K$ in $\R^d$ and define a wedge-type domain by the condition $\{ z \in \Omega: \rho(z) \in K \}$. However, in many cases the study of holomorphic functions on such domains can be reduced to that on the simpler wedges $W(\Omega,E)$. For $\delta > 0$ we also consider a $\delta$-``truncated" wedge 
$$
W_\delta(\Omega,E) =\left\{ z \in \Omega: \rho_j(z) - \delta\sum_{k \neq j} \rho_k < 0, \ j= 1,...,d \right\} \subset W(\Omega,E).
$$
The following result follows by the complexification of a real analytic parametrization of a totally real submanifold.

\begin{prop}
\label{RealAnalytic}
Let $E$ be a real analytic totally real submanifold of dimension $n$ in $\C^n$. For every point $p \in E$ there exists an open neighbourhood $\Omega$ in $\C^n$ and a holomorphic embedding  $\Phi:\Omega \to \C^n$ such that $\Phi(p) = 0$ and $\Phi(E \cap \Omega) = \R^n \cap \Phi(\Omega)$.
\end{prop}

This proposition simplifies many aspects of complex analysis near real analytic totally real submanifolds of maximal dimension. 
If $E$ is merely smooth, then a more subtle result holds: there exists a diffeomorphism $\Phi$ which takes 
$E$ to $\R^n$  and such that $\overline\partial\Phi$ vanishes to infinite order on $E$.

In the study of totally real submanifolds the following {\it gluing disc argument} is often quite helpful. It was introduced in \cite{Pi1} and then used by many authors. Without loss of generality, we may assume that in a neighbourhood $\Omega$ of the origin a smooth totally real manifold $E$ is defined by the equation $x = r(x,y)$, where a smooth vector function $r = (r_1,...,r_n)$ satisfies the conditions $r_j(0) = 0$, $\nabla r_j(0) = 0$. Fix a positive noninteger $s$ and consider for a real function $u \in C^s(\partial\D)$  the Hilbert transform $H: u \to H(u)$. It is uniquely defined by the conditions that  the function $u + iH(u)$ is the trace of a function holomorphic on $\D$ and the integral average of $H(u)$ over the circle is equal to $0$. This is a classical linear singular integral operator; it is  bounded on the space $C^s(\partial\D)$. Let $S^+ = \{ e^{i\theta}: \theta \in [0,\pi] \}$ and $S^- = \{ e^{i\theta}: \theta \in ]\pi, 2 \pi[ \}$ be the semicircles. Fix a $C^\infty$-smooth real function $\psi_j$ on $\partial\D$ such that $\psi_j\vert S^+ = 0$ and  $\psi_j\vert S^- < 0$, $j=1,...,n$. Set $\psi = (\psi_1,...,\psi_n)$. Consider the {\it generalized Bishop equation} 
\begin{eqnarray}
\label{Bishop1}
u(\zeta) = r(u(\zeta),H(u)(\zeta) + c) + t\psi(\zeta), \,\, \zeta \in \partial\D ,
\end{eqnarray}
where $c \in \R^n$ and $t = (t_1,...,t_n)$, $t_j \ge 0$, are real parameters. It follows by the implicit function theorem that this equation admits a unique solution 
$u(c,t) \in C^s(\partial\D)$ depending smoothly on the parameters $(c,t)$. Consider now the analytic discs $f(c,t)(\zeta) = P_\D(u(c,t)(\zeta) + iH(u(c,t))(\zeta))$, where $P_\D$ denotes the Poisson operator of harmonic extension to $\D$. The map $(c,t) \mapsto f(c,t)(0)$ (the centres of discs) is of class $C^s$.  Every disc is attached to $E$ along the upper semicircle. It is easy to see that this family of discs fills the wedge $W_\delta(\Omega,E)$ when $\delta > 0$ and a neighbourhood $\Omega$ of the origin are chosen small enough. Indeed, this is immediate when the function $r$ vanishes identically (i.e., $E = i\R^n$), while the general case follows by a 
small perturbation argument.

This construction of gluing analytic discs is flexible enough and has several applications. As an example we prove a version of {\it the edge-of-the-wedge theorem} following \cite{Ai,Tu1}.

 Consider the generic manifolds $E_j = \{ z \in \Omega: \rho_k(z) = 0, k \neq j, k = 1,...,n \}$ of dimension $n+1$. 
 On the unit circle we consider the open arcs $S_j$, $j=1,...,n$, bounded by the points $\{e^{\frac{2\pi j}{n}i}$, $j =  0,...,n-1\}$. Let $\psi_j$ be $C^\infty$-smooth functions on $\partial\D$ such that $\psi_j \vert S_j < 0$ and $\psi_j \vert (\partial\D \setminus S_j) = 0$, $j=1,...,n$. The equation (\ref{Bishop1}) admits a solution in $C^s(\partial\D)$ smoothly depending on the parameters $(c,t)$ in a neighbourhood of the origin in $\R^{2n}$ (note that $t_j$ are not assumed to be positive here). Every analytic disc from the  family $f(c,t)(\zeta)$ obtained as above has the boundary attached to the union $\cup_j E_j$. Furthermore, their centres $f(c,t)(0)$ fill a neighbourhood of the origin in $\C^n$. Indeed the map $(c,t) \mapsto f(c,t)(0)$ has the maximal rank $2n$ in a neighbourhood of the origin (this is obvious when $r = 0$ and hence remains true under small perturbations).  In combination with the approximation result (Theorem~\ref{BaTr}) we obtain 

\begin{prop}
Let $f$ be a continuous CR function on $\cup_j E_j$. Then $f$ extends holomorphically to a neighbourhood of $E$ in $\C^n$.
\end{prop}

Indeed, by the maximum principle (applied along every analytic disc) the approximating family of holomorphic functions converges in a neighbourhood of the origin.

As a corollary we obtain the  edge-of-the-wedge theorem (for a more general result see \cite{Pi3}). 
Introduce the domains 
$\Omega^+ = \{ z \in \Omega: \rho_j > 0, j=1,...,n \}$ and $\Omega^- = \{ z \in \Omega^-: \rho_j(z) < 0, j=1,...,n \}$.

\begin{cor}
Let $f^+$ and $f^-$ be functions holomorphic on the wedges $\Omega^+$ and $\Omega^-$ respectively and continuous up to the edge $E$. If $f^+$ and $f^-$ coincide on $E$, then they extend to a holomorphic function in a neighbourhood of $E$.
\end{cor}

In combination with Proposition~\ref{RealAnalytic} and the Schwarz Reflection Principle this immediately gives the following simple multidimensional version of this principle:

\begin{prop}\label{CR2}
Let $E$ and $E'$ be real analytic totally real manifolds of dimension $n$ and $N$ in $\C^n$ and $\C^N$ respectively. Suppose that 
$f:W(\Omega,E) \to \C^N$ is a holomorphic mapping continuous on $W_\delta(\Omega,E) \cup E$ for some $\delta > 0$ and such that $f(E) \subset E'$. Then $f$ extends holomorphically to a neighbourhood of $E$.
\end{prop}

A smooth version of this result also holds but requires some additional  technical tools.

\begin{prop}
\label{CR3}
Let $W(\Omega,E)$ be a wedge in $\C^n$ with a $C^\infty$-smooth totally real edge $E$ of dimension $n$. Suppose that $f:W(\Omega,E) \to \C^N$ is a holomorphic mapping such that the cluster set $C_{W(\Omega,E)}(f;E)$ is contained in a 
$C^\infty$-smooth totally real manifold $E'$ of dimension $N$. Then for every $\delta > 0$ 
the mapping $f$ extends to $W_\delta(\Omega;E) \cup E$ as a $C^\infty$-smooth mapping.
\end{prop}
We sketch the proof based on the ideas of Pinchuk-Hasanov \cite{PiHa} (for details see \cite{CoSu1}). The first step is to establish the result for $n=1$, i.e., when $f$ is an analytic disc. This is a combination of  Proposition \ref{CCS} and Chirka's boundary regularity theorem for analytic discs \cite{Ch2}. The second step is to apply the above construction of filling $W(\Omega,E)$ with analytic discs glued to $E$ along the upper semicircle. Since H\"older constants are uniform with respect to the parameters, this implies the H\"older continuity of $f$ up to $E$. In the last step we use the smooth version of Proposition~\ref{RealAnalytic}. Let 
$\Phi$ (resp. $\Psi$) be a (local) diffeomorphism which takes 
$E$ to $\R^n$ (resp. $E'$ to $\R^N$) and such that $\overline\partial\Phi$ and $\overline\partial\Psi$ vanish to infinite order 
on $E$ (resp. on $E'$). We can apply 
the usual Reflection Principle to the mapping $\Psi \circ f \circ \Phi^{-1}$. This gives two functions in the opposite wedges with the edge $\R^n$.
The functions are continuous up to $\R^n$, coincide there, and have the property that the $\overline\partial$-part of their differential vanishes to a suitable order  on $\R^n$. 
But then these functions are $C^\infty$ smooth up to the edge $\R^n$. This is a very special case of the general elliptic regularity of the $\overline\partial$-operator. In our case it can be directly proved by slicing with complex linear discs and using regularity of the Cauchy integral transform $f \mapsto (2\pi i)^{-1} f * (1/\zeta)$ on $\D$.

\section{Geometry of real analytic hypersurfaces}\label{s.geom}

Real hypersurfaces in $\C^n$, $n > 1$, have nontrivial geometry induced by the complex structure of the ambient space. 
This is the main reason for rigidity of holomorphic mappings between domains in $\C^n$. In this section we describe 
classical methods used in the investigation of rigidity properties of holomorphic mappings near boundaries of domains.

\subsection{Complexification, Segre varieties, and differential equations} 
We first introduce an important family of local biholomorphic invariants of a Levi nondegenerate real analytic hypersurface 
$\Gamma$ in $\C^n$, $n >1$. This is a family of complex hypersurfaces 
called {\it Segre varieties} of~$\Gamma$. One can view them as (the graphs of) solutions of a holomorphic second order 
PDE system with a completely integrable prolongation to the space of 1-jets. When $n=2$ such a system becomes a 
second order holomorphic ODE and the Segre family consists of complex curves. The biholomorphic maps of $\Gamma$ 
are precisely Lie symmetries of its Segre family. Thus, the geometry of real analytic hypersurfaces is closely related to the 
geometry of holomorphic ODEs and PDEs. This fundamental correspondence, discovered by Segre \cite{Se1}, inspired 
E.~Cartan \cite{Car} to study the geometry of real hypersurfaces in $\C^2$ in analogy with the geometry of a second 
order ODE developed by the school of S.~Lie \cite{Tr}. The approach of E.~Cartan is very different and is based on his 
equivalence method for Pfaffian systems.

All considerations of this section are local so the results should be understood in terms of the germs of the analytic objects 
involved. To simplify the notation, we will not use the language of germs, so the reader should keep in mind the locality 
assumption.

 Let $\G$ be a real analytic hypersurface in a neighbourhood of $0 \in \G$ in $\C^n$. Then $\G = \{ z: \rho(z,\overline{z}) = 0 \}$, where $\rho$ is a local defining real analytic function. For $w$ close enough to the origin we consider the complex hypersurface
\begin{eqnarray}
\label{Segre1}
Q_w = \{ z: \rho(z,\overline{w}) = 0 \} .
\end{eqnarray}
This hypersurface is called {\it the Segre variety} of the point $w$ (associated with $\G$). The collection of all Segre varieties 
is called {\it the Segre family} of $\G$. More generally, if $\G$ is any real analytic set defined as the zero set of the vector function $\rho$, its Segre varieties are complex analytic subsets of $\C^n$ also defined by (\ref{Segre1}).

The following basic properties of the Segre family can be easily checked:
\begin{itemize}
\item[(i)] $z \in Q_z$ if and only if $z \in \G$.
\item[(ii)] $z \in Q_w$ if and only if $w \in Q_z$.
\item[(iii)] Let $F$ be a holomorphic mapping in a neighbourhood of the origin such that $F(\G) \subset \G'$, where $\G'$ is another real analytic hypersurface. Then $F(Q_w) \subset Q_{F(w)}'$, where $Q_\bullet'$ denotes the Segre family of $\G'$.
\end{itemize}
Property (iii) means that the Segre family is invariant with respect to biholomorphic mappings. In the one-dimensional case Segre varieties are points and so (iii) becomes the classical Schwarz Reflection Principle. In higher dimensions this property leads to far reaching consequences. For applications it is convenient to state (iii) in a more general form.

\begin{lemma}
\label{SegreLemma}
Let $0\in \Gamma$ be a real analytic hypersurface in $\C^n$ and $0\in \G'$ be a real analytic subset of $\C^N$. 
Let $F$ be a holomorphic mapping such that $F(0)=0$ and $F(\G) \subset \G'$. Then $F(Q_w) \subset Q_{F(w)}'$, 
where $Q_\bullet'$ denotes the Segre family of $\G'$.
\end{lemma}

For the proof let $\G' = \{ z' \in \C^N: \phi_j(z',\overline{z}') = 0, j = 1,...,d \}$. Then $\phi_j(F(z),\overline{F(z)}) = 0$ whenever
$\rho(z,\overline{z}) = 0$. Therefore, $\phi_j(F(z),\overline{F(z)}) = \lambda_j(z,\overline{z})\rho(z,\overline{z})$, 
where $\lambda_j$ is a real analytic function in a neighbourhood of the origin. It follows that $\phi_j(F(z),\overline{F(w)}) = \lambda_j(z,\overline{w})\rho(z,\overline{w})$ for $(z,w)$ close to the origin in $\C^n \times \C^N$. This proves the lemma.

We now draw a connection between the complex geometry of real analytic hypersurfaces and the geometry of analytic differential equations and projective connections. The main idea is that the Segre family is a general set of solutions of some second order PDE system (or a single second order ODE when $n=2$).

Let us discuss in some detail the case of dimension 2. We begin with the basic example of 
$\G =\{z_2 + \overline{z}_2 + z_1 \overline{z}_1 = 0\}$, an unbounded realization of the unit sphere in $\mathbb C^2$. 
The Segre family has the form $Q_w= \{ z\in \mathbb C^2 : z_2 + \overline{w}_2 + z_1 \overline{w}_1 = 0 \}$
This is simply the family of all complex lines in $\C^2$ (except the ``vertical" lines $z_1 = $const).
We view every $Q_w$ as the graph of a complex affine function $z_2 = h(z_1)$ that depends on two complex parameters $w_1$ and $w_2$. We treat $z_1$ as an independent variable and $z_2$ as the dependent one. Then the Segre family is the set of graphs of all solutions of the ordinary differential equation $\ddot z_2 = 0$.
In the general case we can assume that $M = \{ z: \Re z_2  = \phi(z_1,\overline{z}_1, \Im z_2) \}$,
where $\nabla \phi(0) = 0$. By the Implicit Function Theorem, $Q_w = \{z : z_2 = h(z_1,\overline{w}_1,\overline{w}_2) \}$
for some holomorphic function $h$. Again we view $z_2$ as the dependent variable and $z_1$ as the independent one. Applying  the Chain Rule we obtain
$d^jz_2/dz_1^j = (\partial^j h/\partial z_1^j)(z_1,\overline{w}_1,\overline{w}_2)$, $j=1,2$.
By the Implicit Function Theorem  we represent  the parameters $(w_1,w_2)$ as the functions of $(z_1, z_2,\dot z_2)$ and  obtain the holomorphic ODE
\begin{eqnarray}
\label{Segre7}
\ddot z_2 = F\left(z_1,z_2,\dot z_2\right) .
\end{eqnarray}
The Segre family of the real hypersurface $\G$  is precisely the set of the graphs of the solutions of~(\ref{Segre7}). 

The invariance property (iii) of the Segre family means that a biholomorphism $f$ of $\G$ sends the graph of a solution of (\ref{Segre7}) to the graph of another solution. But this means precisely that $f$ is a (point) Lie symmetry of the equation (\ref{Segre7}).  Therefore, the classical theory of Lie symmetries can be applied in order to study biholomorphisms of real analytic hypersurfaces.

As an example, consider a holomorphic differential equation
\begin{eqnarray}
\label{Segre8}
(S): \ddot u = F(x,u,\dot u) ,
\end{eqnarray}
where $x$ denotes an independent complex variable and $F$ is a holomorphic function. A   {\it symmetry group} $Sym(S)$ of the equation $(S)$ is a (maximal) local (Lie) group $G$ acting on a domain in~$\C^2$ such that the following holds: for every solution $u(x)$ of $(S)$ and every $g \in G$ the image (if defined)  of the graph of $u$ by $g$ is the graph of some solution of $(S)$ which we denote by $g_*u$. A holomorphic vector field 
  \begin{eqnarray}
  \label{Segre9}
  X = \theta\frac{\partial}{\partial x} + \eta\frac{\partial}{\partial u}
  \end{eqnarray}
is called an {\it infinitesimal Lie symmetry} of $(S)$ if it belongs to the Lie algebra of $Sym(S)$, i.e., generates a one-parameter group of point Lie symmetries of $(S)$. Denote by $j_x^m(u)$ the $m$-jet of $u$ at $x$. In the most important case $m=2$ we
set $u_1 = u_x$ and $u_{11} = u_{xx}$. Then  $j^2_x(u) = (x,u,u_1,u_{11})$, and so $(x,u,u_1,u_{11})$ are natural 
coordinates on this jet space.

Every point Lie symmetry $g$ canonically extends to $J^m(1,1)$ as a biholomorphic mapping $g^{(m)}$ defined as follows: $g^{(m)}$ associates to $J_x^m(u)$ the jet $j_{u(x)}^m(g_*u)$. In particular, a one-parameter group of symmetries generated by a vector field $X$  lifts to $J^m(1,1)$. A vector field $X^{(m)}$ on $J^m(1,1)$ which generates this lift is called the  prolongation of order $m$ of $X$. The classical Lie theory gives powerful  tools for the study of Lie symmetries which are particularly convenient in the infinitesimal case. For $m=2$ we have
\begin{eqnarray}
  \label{Segre10}
  X^{(2)} = X + \eta_1\frac{\partial}{\partial u_1} + \eta_{11}\frac{\partial}{\partial u_{11}} .
 \end{eqnarray}
The general Lie theory gives the following expressions for the coefficients:
\begin{eqnarray}
\label{Segre11}
\eta_1 = \eta_x + (\eta_u - \theta_x)u_1 - \theta_u(u_1)^2 ,
\end{eqnarray}
\begin{eqnarray}
\label{Segre12}
\eta_{11} = \eta_{xx} + (2\eta_{xu} - \theta_{xx})u_1 + (\eta_{uu} - 2\theta_{xu})(u_1)^2 - \theta_{uu}(u_1)^3 + (\eta_u - 2\theta_x)u_{11} - 3\theta_u u_1 u_{11} .
\end{eqnarray}
Equation (S) of the form (\ref{Segre8}) defines a complex hypersurface $(S_2)$ in $J^2(1,1)$ by the equation $u_{11} = F(x,u,u_1)$. The fundamental principle of the Lie theory states that $X$ is an infinitesimal symmetry of $(S)$ if and only if the vector field $X^{(2)}$ is tangent to $(S_2)$, that is, 
\begin{eqnarray}
\label{Segre13}
X^{(2)}(u_{11} - F(x,u,u_1)) = 0 \,\,\, \mbox{for} \,\,\,(x,u,u_1,u_{11}) \in (S_2).
\end{eqnarray}

Consider the expansion $F(x,u,u_1) = \sum_{\nu \ge 0} f_{\nu}(x,u)(u_1)^\nu$. Plugging it into (\ref{Segre13}) and comparing the coefficients of the powers of $u_1$ we obtain a system of PDEs of the form 
\begin{eqnarray*}
L D^2(\theta,\eta) = G(x,u,D^1(\theta,\eta)) .
\end{eqnarray*}
Here $D^j$ denotes the set of the partial derivatives of the map $(\theta,\eta)$ of order $j$, $G$ is an analytic function, and 
$L$ is a matrix  with constant coefficients. Applying to this system the partial derivatives in $x$ and $u$, we obtain after a direct computation that 
\begin{eqnarray}
\label{Segre14}
D^3(\theta,\eta) = H(x,u,D^1(\theta,\eta),D^2(\theta,\eta)) ,
\end{eqnarray}
for some analytic function $H$. This implies that every infinitesimal Lie symmetry of $(S)$ is determined by its second order jet at a given point. In particular, $\dim Sym(S) \le 8$.

Consider again the equation (\ref{Segre7}) describing the Segre family of $\G$. The group of local biholomorphisms of $\G$ is embedded into the symmetry group of (\ref{Segre7}) as a totally real subgroup of maximal dimension (see more details in \cite{Su3}). As a consequence we obtain that the dimension of the real Lie group of biholomorphisms of $\G$ is bounded above by 8.

In higher dimensions ($n > 2$) the Segre family of a Levi nondegenerate real analytic hypersurface $\Gamma$ is described by a PDE system 
\begin{eqnarray}
\label{Lie1}
u_{x_ix_j} = F_{ij}(x,u,u_x) ,
\end{eqnarray}
where $x \in \C^{n-1}$ and $u_x$ denotes the set of the first order partial derivatives of the function $u = u(x)$.
This system is completely integrable: its lift to the first order space of jets is a 
first order PDE system which satisfies  the Frobenius integrability conditions. For more details we refer the reader to papers
\cite{Me3, Su3}, which are devoted to the study of Lie symmetries of PDE systems and Segre families. 

\subsection{Equivalence problem for real hypersurfaces, I: Moser's Approach}

Let $\G$ be a real analytic strictly pseudoconvex hypersurface in $\C^n$ containing the origin. 
All consideration will be local. We use the notation $z = ('z,z_n)$, $'z \in \C^{n-1}$. By the Implicit Function Theorem (after a permutation of coordinates),  $\G$ can be written as the graph
\begin{eqnarray}
\label{Moser1}
y_n =  F('z,'\overline{z},x_n)
\end{eqnarray}
of a  real analytic function $F$.  Given a point $p \in \G$ consider a real analytic curve $\gamma$ in $\G$ with a  
parametrization $z = z(\tau)$,  $\tau \in (-\tau_0,\tau_0)$. Assume that $\gamma$ passes through  $p$ in a noncomplex 
tangential direction, i.e., $z(0) = p$ and the vector $\dot z(0)$ is not in $H_p\G$. 

In a neighbourhood of $p$ there exists a biholomorphic map $z^* = h(z)$ taking the curve $\gamma$ to the real interval  $'z = 0$, $z_n = \tau$ (we drop the asterisk for simplicity), that is, $h(z(\tau)) = ('0,\tau)$. Furthermore, in the new coordinates $\G$ is given by the equation
\begin{eqnarray}
\label{Moser3}
y_n = \vert 'z \vert^2 + \sum_{k,l \ge 2} F_{kl}('z,'\overline{z},x_n),
\end{eqnarray}
where $F_{kl}$ are real homogeneous polynomials of degree $k$ in $'z$ and $l$ in $'\overline{z}$ with coefficients analytic in $x_n$. 

Of course, such a change of coordinates  is not unique. First of all this is due to the freedom in the choice of the curve $\gamma$ and its parametrization. Moreover, consider a transformation
\begin{eqnarray}
\label{Moser2}
('z,z_n) \mapsto (U(z_n)'z,z_n) .
\end{eqnarray}
Here $z_n \mapsto U(z_n)$ is a $(n-1)\times (n-1)$ nondegenerate holomorphic  matrix function in $z_n$ which is unitary for $z_n = \tau$. This transformation fixes $\gamma$ and preserves (\ref{Moser3}).

Using tensor notation we write
$$F_{kl} = \sum_{1 \le \alpha_\nu,\beta_\mu \le n-1} a_{\alpha_1...\alpha_k\beta_1...\beta_l}
z_{\alpha_1}...z_{\alpha_k}\overline{z}_{\beta_1}...\overline{z}_{\beta_l} .$$
Here we assume that the coefficients $a_{\alpha\beta}$ do not change under permutations of indices $\alpha_\nu$ and $\beta_\mu$. For $k,l \ge 1$, we put
$$tr F_{kl} = \sum b_{\alpha_1...\alpha_{k-1}\beta_1...\beta_{l-1}}z_{\alpha_1}...z_{\alpha_{k-1}}
\overline{z}_{\beta_1}...\overline{z}_{\beta_{l-1}} ,$$
with
$$ 
b_{\alpha_1...\alpha_{k-1}\beta_1...\beta_{l-1}} = \sum_{\alpha_k = \alpha_l} a_{\alpha_1...\alpha_k\beta_1...\beta_l} .
$$
Moser \cite{CM} proved that after a biholomorphic change of coordinates one can additionally achieve in (\ref{Moser3}) the conditions
\begin{eqnarray}
\label{Moser4}
tr F_{22} = (tr)^2 F_{32} = (tr)^3 F_{33} = 0 .
\end{eqnarray}
Representation (\ref{Moser3}), (\ref{Moser4}) is called {\it the (Moser) normal form} of $\G$. A real analytic curve $\gamma$, 
which in the  normal form has the equation $'z = 0$, $y_n = 0$, is called a {\it chain}.

Conditions (\ref{Moser4}) may be viewed geometrically. First we note that for a given point $p \in \Gamma$ there exists a unique 
chain passing trough $p$ in a prescribed noncomplex tangential direction. Furthermore, if $\Gamma$ is given by (\ref{Moser3}), 
then the line $'z = 0$, $y_n = 0$ is a chain iff $(tr)^2 F_{32} = 0$.  Let now $A$ be a unitary $(n-1) \times (n-1)$ matrix. There exists a unique mapping (\ref{Moser3}) such that $U(0) = A$ and in the new coordinates $tr F_{22} = 0$. This matrix $A$ can be viewed as a new choice of an orthonormal basis in $H_0\Gamma$. Finally, consider an admissible reparametrization $('z,z_n) \mapsto (\sqrt{\dot q(z_n)}'z,q(z_n))$ of~$\gamma$. Here $q(0) = 0$ and $q(\overline{z}_n) = \overline{q(z_n)}$.
Such a change of coordinates preserves (\ref{Moser3}) and the conditions $tr F_{22} = 0$, $(tr)^2 F_{32} = 0$. 
If $q$ additionally satisfies a certain third order ODE then the condition $(tr)^3F_{33} = 0$ also holds. 
Since $q(0) = 0$, the solution to this equation is uniquely determined by its first and second order derivatives at the origin.

Thus, the normalization of $\Gamma$ depends on the following data:
\begin{itemize}
\item[(i)] the point $p \in \Gamma$ corresponding to the origin in the normal form;
\item[(ii)] a noncomplex tangential direction at $p$ defining the chain $\gamma$ which has the equation $'z = 0$, $y_n=0$ in the normal form;
\item[(iii)] the choice of an orthonormal basis in $H_p\G$;
\item[(iv)] two real parameters fixing the parametrization of $\gamma$.
\end{itemize}

The main result of Moser's theory can be stated as follows.
\begin{theorem}
For each choice of the initial data (i)-(iv) there exists a unique biholomorphic mapping $h$ taking $\G$ to the normal form. 
\end{theorem}
The normal form of the sphere is $y_n = \vert 'z \vert^2$. Every choice of the initial condition (i)-(iv) determines a unique linear-fractional automorphism of the sphere. Therefore, by Moser's theorem every local automorphism is global. We obtain the Poincar\'e-Alexander theorem, see the next section for a detailed discussion.

Many applications of Moser's theory are contained in the expository papers by Vitushkin \cite{Vi,Vi2}.

\subsection{Equivalence problem, II: the Cartan-Chern Approach} Cartan's equivalence problem (in local form)  
can be stated as follows. Let $U$ and $V$ be open subsets in $\R^n$. 
Suppose that $\omega_U = (\omega_U^1,...,\omega_U^n)$ and $\Omega_V = (\Omega_V^1,...,\Omega_V^n)$ 
are co-frames (bases of 1-forms) on $U$ and $V$ respectively. Consider a prescribed linear group $G$. The problem 
is to determine all diffeomorphisms 
$$f: U \to V,$$
satisfying
$$f^*\Omega_V = \omega_U\gamma_{UV} \,\,\,\mbox{with}\,\, \gamma_{UV} \in G.$$
In the original work of E.~Cartan, the group $G$ is allowed to vary from point to point. Note that many natural equivalence 
problems in differential geometry (for Riemannian structures, differential equations, CR structures, etc.) can be represented 
in this form for an appropriate choice of co-frames and the group $G$.

The approach of E.~Cartan to the solution of the above equivalence problem is based on the observation that the problem has a (relatively) simple solution in the case of the trivial group $G = \{ e \}$. Even then a complete solution of the equivalence problem is not quite explicit. In most cases this method provides a number of geometric invariants of the problem at hand (a common application of these invariants is to show that the equivalence problem does not admit any solution). The equivalence problem is considered to be solved if it is reduced to the problem with $G = \{ e \}$. The main idea of this reduction consists of introduction of a (finite) sequence of larger spaces reducing the group $G$ at every step.

The procedure proposed by E.~Cartan comprises several iterated steps. First, the co-frames $\omega$ and $\Omega$ must be extended in an equivariant way to $U \times G$ and $V \times G$. There,  the first structure equations for $d\omega$ and $d\Omega$ can be written. These equations contain the so-called torsion terms and the special procedure of absorption of these terms allows one to reduce the group. Note that this requires the expansion of the initial system to a higher bundle.

In the case of Levi nondegenerate hypersurfaces in $\C^n$ (not necessarily real analytic but sufficiently smooth) this approach leads to a complete solution of the local equivalence problem. This was achieved by E.~Cartan \cite{Car} for $n=2$ and by Chern~\cite{CM}
and Tanaka~\cite{Ta} in all dimensions. The adjacent equivalence problem for  systems (\ref{Lie1}) was solved by Hachtroudi \cite{Hach} using Cartan's method. This problem can be also viewed in terms of Cartan's theory of the projective connections. This approach to the Segre geometry is developed by Chern \cite{C} and Burns-Shnider \cite{BuSh3}.

Each of the methods described in this section gives useful information concerning local properties of biholomorphic mappings between strictly pseudoconvex or Levi nondegenerate hypersurfaces: precise bounds on the dimension of the automorphism group, efficient parametrization of biholomorphisms by their second order jets, etc. Notice that the methods of Moser and Cartan-Chern admit some generalizations to a wider class of real hypersurfaces (with special Levi degeneracies) or to real manifolds of higher codimension. However, in these cases the approach based on Segre varieties often turns out to be the most convenient and flexible. For example, the approach based on the dynamical description of the Segre family by an analytic ODE was  recently extended to a wider class of hypersurfaces (with  the degenerate Levi form) by Kossovskiy-Shafikov \cite{KoSha1,KoSha2}. Using the technique of the local theory of analytic dynamical systems and meromorphic differential equations they studied the geometry of Levi degenerate hypersurfaces and their formal (i.e., with respect to maps given by formal power series) and analytic CR transformations. In particular, they proved in \cite{KoSha2} that formal  and holomorphic equivalence of real analytic hypersurfaces in $\C^n$ do not coincide (for Levi nondegenerate hypersurfaces they coincide by Moser's theory). This is a consequence of the classical phenomenon of the local 
theory of analytic dynamical systems where the formal and analytic classifications are different. Further results in this direction are  obtained recently by Kossovskiy-Lamel \cite{KoLa}.

\section{Holomorphic mappings of strictly pseudoconvex domains and scaling}
\label{s.refl}

We begin with the special case of the unit ball $\B^n \subset \C^n$, $n > 1$, and the Poincar\'e -Alexander rigidity 
phenomenon~\cite{Po, Al1}:
let $U$, $V \subset \C^n$ be neighbourhoods of points $p, q \in \partial \B^n$ respectively, and $f: U \to V$, $f(p)=q$, be a biholomorphic (or even nonconstant holomorphic) map which takes $U \cap \partial\B^n$ to $V \cap \partial \B^n$; then $f$ extends to an automorphism of $\B^n$.

This result was used by Alexander \cite{Al2} and Rudin \cite{Ru1} to prove the following 

\begin{theorem}
\label{AlRu}
Any proper holomorphic self-map of the unit ball $\B^n$, $n>1$, is an automorphism of~$\B^n$.
\end{theorem}

The Poincar\'e -Alexander phenomenon and Theorem \ref{AlRu} illustrate a great difference between biholomorphic and proper holomorphic mappings in one and several complex variables. The following result  obtained by C.~Fefferman \cite{Fe} has been influential for the development of the theory of holomorphic mappings.

\begin{theorem}\label{FefTheo}
Let $f: \Omega \to \Omega'$ be a biholomorphic mapping between two strictly pseudoconvex domains in $\C^n$ with $C^\infty$-smooth boundaries. Then $f$ extends to $\overline\Omega$ as a mapping of class $C^\infty$.
\end{theorem}

The original proof of this theorem was long and complicated, but it stimulated intensive research on the boundary regularity of proper holomorphic mappings which led to the discovery of new methods and results. In this section we present the main steps of a different,
more elementary approach to the proof of Theorem \ref{FefTheo}.

\subsection{The Schwarz Reflection Principle. } 
The following theorem may be considered as the Schwarz Reflection Principle in $\C^n$, $n>1$.

\begin{theorem}\label{SchwarzTheo1}
Let $\Omega$ be a one-sided neighbourhood of a strictly pseudoconvex real analytic hypersurface $\G$ in $\C^n$. Let also $\G'$ be a real analytic strictly pseudoconvex hypersurface in $\C^n$. Suppose that $f: \Omega \to \C^n$ is a holomorphic mapping of class $C^1(\Omega \cup \G)$ such that $f(\G) \subset \G'$. Then $f$ extends holomorphically to a full neighbourhood of $\G$ in $\C^n$.
\end{theorem}

\begin{cor}
Let $\Omega$, $\Omega'$ be strictly pseudoconvex domains with real analytic boundaries and $f:\Omega \to \Omega'$ be a proper holomorphic map which extends to $\overline\Omega$ as a $C^1$ map. Then $f$ extends holomorphically to a neighbourhood of $\overline\Omega$.
\end{cor}

We present here two different proofs of this theorem, both of which are  local. The first one was obtained in \cite{Pi1} and \cite{Lewy}.

Fix a point $p \in \G$.
Let $\rho$ and $\psi$ be strictly plurisubharmonic real analytic local defining functions of $\G$ and $\G'$ near $p$ and $f(p)$ respectively. As usual, we consider them as power series in $z$ and $\overline{z}$. Then 
\begin{eqnarray}
\label{Schwarz1}
\psi(f(z),\overline{f(z)}) = 0,
\end{eqnarray}
whenever 
\begin{eqnarray}
\label{Schwarz2}
\rho(z,\overline{z}) = 0.
\end{eqnarray}
Let $X_1,...,X_{n-1}$ be a basis in the space of tangential Cauchy-Riemann operators on $\G$.
Applying them to (\ref{Schwarz1}) we obtain via the Chain Rule,
\begin{eqnarray}
\label{Schwarz3}
X_j \psi(f(z),\overline{f(z)}) = \sum_{k=1}^{n} \frac{\partial \psi}{\partial w_k}(f(z),\overline{f(z)}) X_jf(z) = 0 , \,\, j = 1,...,n-1.
\end{eqnarray}
It suffices to consider the case when $f$ is not constant. Since both hypersurfaces are strictly pseudoconvex we conclude that the tangent mapping $df(p)$ is nondegenerate. Together with  the nondegeneracy of the Levi form of $\psi$ at $f(p)$, this allows us to apply the Implicit Function Theorem to the system (\ref{Schwarz1}), (\ref{Schwarz3}). We obtain
\begin{eqnarray}
\label{Schwarz4}
\overline{f(z)} = H(f(z),X_1f(z),...,X_{n-1}f(z)), \,\,\,z \in \G,
\end{eqnarray}
where $H$ is a holomorphic function in all variables. Let $\C \ni \zeta \mapsto l_c(\zeta)$ be a family of parallel complex affine lines depending on a parameter  $c \in \C^{n-1}$ that are transverse to $\G$ near $p$. Then the intersection of every line with $\G$ is a real analytic curve $\gamma_c$ in $\C$. Every coefficient of $X_j$ is a real analytic function and its restriction on $\gamma_c$ extends to a function holomorphic in a neighbourhood (whose size is independent of $c$) of $\gamma_c$. Replacing in (\ref{Schwarz4}) these coefficients with such holomorphic extensions, we obtain by the one-dimensional Schwarz Reflection Principle that the restriction of $f$ on every linear section $\Omega \cap l_c$ extends holomorphically past $\G$. Then the result follows by the Hartogs Lemma.

The second proof, due to Webster \cite{We3}, is shorter. Since the hypersurfaces are strictly pseudoconvex,  the projectivizations of their holomorphic tangent bundles $\PP H(\G)$ and $\PP H(\G')$ are totally real submanifolds of maximal dimension $2n-1$. Then the mapping $z \mapsto (z,df(z))$ is holomorphic on a wedge type domain with the edge $\PP H(\G)$, is continuous up to the edge, and takes it to $\PP H(\G')$. Hence, this mapping extends holomorphically to a neighbourhood of the edge by the Reflection Principle
(see Proposition~\ref{CR2}).

These two proofs of Theorem \ref{SchwarzTheo1} represent two different types of the Reflection Principle: analytic and geometric. 
They both require additional assumptions on the regularity of $f$ on the hypersurface $\Gamma$. The analytic approach needs 
$C^1$ regularity, while the geometric one requires slightly weaker regularity, namely continuity of the lift of $f$ up to $\PP H(\partial\Omega)$. Furthermore, the geometric Reflection Principle admits smooth generalizations. A $C^\infty$-smooth version 
of Theorem~\ref{SchwarzTheo1} was obtained by Nirenberg-Webster-Yang \cite{NiWeYa}: if the hypersurfaces $\G$ and $\G'$ are merely $C^\infty$-smooth, then $f$ necessarily is of class $C^\infty(\Omega \cup \G)$. The proof follows along the same lines 
with the application of Proposition \ref{CR3}.

\subsection{Continuous extension} We need the following version of the classical Hopf lemma.

\begin{prop}
\label{Hopf}
Let $\Omega$ be a bounded domain with $C^2$-smooth boundary in $\C^n$ and let $K$ be a compact subset of $\Omega$.
For every constant $L > 0$ there exists a constant $C = C(K,L) > 0$ with the following property: if a function  $u \in PSH(\Omega)$ is such that $u(z) < 0$ for every $z \in \Omega$ and $u(z) \le - L$ for all $z \in K$, then $\vert u(z) \vert \ge C dist(z,\partial\Omega)$ for each $z \in \Omega$.
\end{prop}

One of the first results on the boundary behaviour of holomorphic mappings (see \cite{Kh,Pi2}) is the following 

\begin{theorem}
\label{ExtTheo1}
Let $f: \Omega \to \Omega'$ be a proper holomorphic mapping between two strictly pseudoconvex domains in $\C^n$. 
Then $f$ extends to $\overline\Omega$ as a 1/2-H\"older-continuous mapping.
\end{theorem}

Let $\rho$ and $\psi$ be strictly plurisubharmonic global defining functions of $\Omega$ and $\Omega'$ respectively. The functions $v(z) = \psi(f(z))$ and  $u(p) = \sup\{ \rho(z): f(z) = p \}$ are plurisubharmonic in $\Omega$ and $\Omega'$ respectively. Applying the Hopf Lemma to these functions we obtain that $C\rho(z)   \le \psi(f(z))  \le C^{-1}\rho(z)$ for some constant $C > 0$. This is equivalent to the boundary distance preserving property 
\begin{eqnarray}
\label{Ext1}
C {\,\rm dist}(z,\partial \Omega )   \le {\,\rm dist} (f(z), \partial \Omega')  \le C^{-1} {\,\rm dist}(z,\partial\Omega) .
\end{eqnarray}
From the decreasing property of the Kobayashi-Royden pseudometric and the estimates of  this metric  from above and below, 
it follows that
\begin{eqnarray*}
C \vert df(z)V\vert {\,\rm dist}(f(z),\partial \Omega')^{-1/2} \le      
F_{\Omega'}(f(z),df(z)V) \le F_{\Omega}(z,V) \le C^{-1}\vert V \vert {\,\rm dist}(z,\partial \Omega)^{-1} ,
\end{eqnarray*}
for every point $z \in \Omega$ and every tangent vector $V$. In view of (\ref{Ext1}) this implies the estimate
\begin{eqnarray}
\label{Ext2}
\parallel df(z) \parallel \le C {\rm dist}(z,\partial \Omega)^{-1/2}
\end{eqnarray}
for the operator norm of the differential. The theorem now follows by the classical integration argument of Hardy-Littlewood.

The same proof works with minor modifications when the domain $\Omega$ is merely pseudoconvex: instead of the defining function
$\rho$ one can use the bounded exhaustion function of Diederich-Fornaess (Theorem~\ref{t.df}). The assumptions on $f$ and $\Omega'$ also can be weakened, see Proposition \ref{ExtTheo2}. It was the idea of Diederich-Fornaess \cite{DiFo3} to utilize the Kobayashi-Royden metric instead of the previously used Carath\'eodory metric.

Theorem \ref{ExtTheo1} does not allow immediately to deduce Fefferman's Theorem~\ref{FefTheo} from Theorem~\ref{SchwarzTheo1} and its smooth counterpart. However, it was used by Nirenberg-Webster-Yang \cite{NiWeYa} 
to prove continuity of the lift of $f$ up to $\PP H(\partial\Omega)$. This was done with rather tricky and subtle arguments 
involving  the Julia-Carath\'eodory lemma. The argument was later simplified by Forstneri\v c \cite{For5}. In the next 
subsection we present a more transparent proof using scaling.

\subsection{The scaling method} Let $\Omega$ be a domain with a strictly pseudoconvex boundary of class $C^2$ and a defining function $\rho$ near a point $w^0 \in \partial\Omega$. There exists a neighbourhood $U$ of $w^0$ in $\C^n$ and a family of biholomorphic mapping $h_w:\C^n \to \C^n$, continuously depending on $w \in \partial\Omega \cap U$, so that the following conditions are satisfied:
\begin{itemize}
\item[(i)] $h_w(w) = 0$.
\item[(ii)] The defining function $\rho_w:= \rho \circ h_w^{-1}$ for the domain $h_w(\Omega)$ has the form
$$\rho_w = 2\Re z_n + 2\Re Q_w(z) + H_w(z) + R_w(z) ,$$
where $R_w(z) = o(\vert z \vert^2)$, $Q_w(z) = \sum_{\mu,\nu=1}^n q_{\mu\nu}(w)z^{\mu}z^{\nu}$ and 
$H_w(z) = \sum_{\mu,\nu=1}^n h_{\mu,\nu}(w)z^{\mu}\overline{z}^{\nu}$. Furthermore, $Q_w(z) = 0$ and $H_w(z) = 0$ when $z_n = 0$.
\item[(iii)] Each mapping $h_w$ sends the real normal of $\partial\Omega$ at the point $w$ to the real normal $\{ z_1=...=z_{n-1} = \Im z_n = 0 \}$ of $\partial h_w(\Omega_w)$ at the origin.
\end{itemize}

In applications of this construction usually one can assume that $w^0=0$ and $\partial\Omega$ is already normalized near the origin; therefore, one can assume additionally that $h_{w^0}$ is the identity mapping.

As before, $'z = (z_1,...,z_{n-1})$ so that $z = ('z,z_{n-1})$. Consider a sequence of points 
$\{q^k\}$ in $\Omega$ converging to a point $q \in \partial \Omega$. Denote by $w^k \in \partial \Omega$ the point closest to $q_k$. Set $h^k:= h_{w^k}$  and $\rho_k := \rho_{w^k}$. Set $\delta_k = {\rm dist}(h^k(q^k),\partial h^k(\Omega^k))$. Then $h^k(q^k) = ('0,-\delta_k)$. Consider the dilations 
$$
d^k:('z,z_n) \mapsto (\delta_k^{-1/2} {'z},\delta_k^{-1}z_n) .
$$ 
Finally, define the biholomorphic mappings $D^k:= d^k \circ h^k$. Note that this sequence of biholomorphic mappings is determined by $\Omega$ and the sequence $\{q^k\}$. We call the sequence $\{D^k\}$ {\it the scaling along a sequence $\{q^k\}$}. Let $\Omega_k = D^k(\Omega) = \{ \delta_k^{-1}\rho_k \circ d_k^{-1} < 0 \}$. It is easy to see that the sequence of functions 
$\{\delta_k^{-1}\rho_k \circ d_k^{-1}\}$ converges uniformly on compact subsets of $\C^n$ to the function $2\Re z_n + \vert z'\vert^2$ which defines the domain $\HH$ given by~\eqref{e.H}. As a consequence, the sequence of domains $\{\Omega_k\}$ 
converges to $\HH$ with respect to  the Hausdorff distance.

Scaling along a sequence has many applications. As an example, we conclude the sketch of the proof of Fefferman's mapping theorem using the arguments from \cite{PiHa}. It suffices to show that in the hypothesis of the theorem the lift $(z,df(z))$ 
of $f$ to the tangent bundle extends continuously to $\PP H(\partial \Omega)$. Arguing by contradiction, assume that there exists a sequence of points $\{p^k\}$ in $\Omega$ converging to a boundary point $p$ such that their images $\{q^k\}$ converge to some point $q \in  \partial\Omega'$, but 
the sequence $\{p^k,df(p^k)\}$ does not converge to $\PP H(\partial\Omega')$. Let $\{G^k\}$ and $\{D^k\}$ be the scaling sequences along the sequences $\{p^k\}$ and $\{q^k\}$ respectively. Then one can show that the sequence 
$\{f^k = D^k \circ f \circ (G^k)^{-1}\}$ converges to a holomorphic mapping which is degenerate at some point. On the other 
hand, it is easy to see by the standard normal family argument that the limit map is a biholomorphism of the unit ball. This contradiction proves the theorem.

The idea of using almost holomorphic functions (i.e., functions with asymptotically vanishing $\overline\partial$-operator) and the Reflection Principle in Fefferman's theorem is due to Nirenberg-Webster-Yang \cite{NiWeYa}. It was also used  for hypersurfaces of 
class $C^m$ with noninteger $m> 2$ in \cite{PiHa} which proves that a biholomorphic mapping $f:\Omega \to \Omega'$ between strictly pseudoconvex domains with boundaries of class $C^m$ extends to $\overline\Omega$ as a map of class $C^{m-1}(\overline\Omega)$. A similar result but by different methods was proved by Lempert \cite{Lem1}. Later Khurumov \cite{Hur} proved that, in fact, $f \in C^{m-1/2}(\overline\Omega)$. This result is sharp.

\subsection{Proper and locally proper mappings}
The case when $f:\Omega \to \Omega'$ is a proper mapping can be reduced to the biholomorphic case 
via  the following generalization of Alexander's theorem obtained in \cite{Pi6}.

\begin{theorem}
\label{AlPinTheo}
Let $f:\Omega \to \Omega'$ be a proper holomorphic mapping between two strictly pseudoconvex domains with $C^2$-smooth boundaries. Then $df(z)$ is nondegenerate at every point $z \in \Omega$, i.e., $f$ is locally biholomorphic.
\end{theorem}

Arguing by contradiction, assume that the Jacobian determinant of $f$ vanishes on a complex hypersurface $H$ in $\Omega$. Let $\{p^k\}$ be a sequence of points in $H$ converging to a boundary point and set $q^k = f(p^k)$. Consider the scalings $\{G^k\}$ 
and $\{D^k\}$ along these sequences. Then the sequence $\{f^k = D^k \circ f \circ (G^k)^{-1}\}$ converges to a proper holomorphic mapping $F: \HH \to \HH$ which is a biholomorphism by Alexander's theorem. On the other hand, it follows by the choice of $\{p^k\}$ that $dF$ vanishes at some interior point, which is a contradiction.

\begin{cor}
\label{AlPinCor}
A proper holomorphic self-mapping of a strictly pseudoconvex domain with $C^2$-smooth boundary is a biholomorphism.
\end{cor} 
 
A similar idea allows one to establish the following rigidity phenomenon for CR mappings \cite{PiTs}.

\begin{theorem}
Let $\G$ and $\G'$ be strictly pseudoconvex hypersurfaces in $\C^n$. Suppose that $U$ is a neighbourhood of a point $p \in \G$ and $f: \G \cap U \to \G'$ is a continuous nonconstant CR mapping. Then there exist neighbourhoods  $V$ and $V'$  of  $p$ and $f(p)$ respectively such that 
$f: \G \cap V \to \G' \cap V'$ is a homeomorphism.
\end{theorem} 

We can use the fact that $f$ extends holomorphically to the pseudoconvex one-sided neighbourhood of $p$. The difficulty is that a 
priori $f$ may be not a proper mapping there. The main idea of the proof is to show that the set $f^{-1}(f(p))$ is finite in a neighbourhood of $p$ on $\partial\Omega$. This will imply that $f$ is proper and will reduce the problem to the previous theorem. Arguing by contradiction suppose that this set contains a sequence converging to $p$ and apply again the scaling (see \cite{CoPiSu1} 
for a more general scaling result needed here). Then one can show that the limit map is an automorphism of the ball and at the same time is degenerate at some point, which is a contradiction.

\section{Extension of germs of holomorphic mappings}

In this section we present the results which generalize and develop the rigidity phenomenon discovered by Poincar\'e and Alexander. 

A real hypersurface $\Gamma = \rho^{-1}(0)$ in $\C^n$ is called {\it real algebraic} if it is defined by a real polynomial~$\rho$. 
For an open set $U\subset \mathbb C^n$, a holomorphic map $F: U \to \mathbb C^n$ is called algebraic if its graph is contained in 
an algebraic subvariety of $\mathbb C^{2n}$ of dimension $n$.
The following result is due to Webster \cite{We1}:

\begin{theorem}
\label{GermsTheo1}
Let $\G$ and $\G'$ be Levi-nondegenerate real algebraic hypersurfaces in $\C^n$ of degree $m$ and $m'$ respectively, and let 
$U$ be a neighbourhood in $\C^n$ of a point $p \in \G$. Suppose that $f:U \to \C^n$ is a holomorphic mapping with a 
nondegenerate differential at $p$ and such that $f(\G) \subset \G'$.
Then $f$ extends to $\C^n$ as an algebraic mapping of degree bounded above by a constant depending only on $n$, $m$, and $m'$.
\end{theorem}

Applying complex conjugation, we rewrite the equation (\ref{Schwarz4}) in the form $f(z) = G(z,\overline{z}, \overline{F(z)})$, where $F = (f,J_f)$. Here $G$ is an (holomorphic) algebraic function and $J_f$ denotes the Jacobian matrix of $f$ viewed as a $\C^{n^2}$-valued map.  It follows that the mapping $F$ takes $\G$ to the real algebraic set $M = \{ (z,\zeta,\omega) \in \C^n \times \C^n \times \C^{n^2}: \zeta = G(z,\overline{z},\overline{\omega})\}$. By Lemma \ref{SegreLemma}, $F(Q_w)$ is contained in $Q_{F(w)}'$ (the Segre variety defined by $M$), which implies that the restriction of $f$ to $Q_w$ is an algebraic map (of controlled degree). 
Consider $n$ families of transverse Segre varieties for $\Gamma$. After a local biholomorphic and algebraic change of coordinates one can transform them to families of parallel coordinate hyperplanes. Now the classical theorem on separate algebraicity \cite{BM} can be applied (see \cite{SharSu} for details).

We begin the discussion of the analytic case with the result of \cite{Pi5}.

\begin{theorem}\label{sphereTheo}
Let $\Gamma$ be a (connected) real analytic strictly pseudoconvex hypersurface in $\C^n$, $n > 1$, $U$ a neighbourhood of a 
point $p \in \Gamma$, $f: U \longrightarrow \C^n$ a nonconstant holomorphic mapping and assume that 
$f(U \cap \Gamma) \subset \partial\B^n$. Then $f$ can be continued along any path on $\Gamma$ starting at $p$ as a locally biholomorphic mapping.
\end{theorem}

\begin{cor}
\label{ExtSpherical}
Let $\Omega$ be a bounded strictly pseudoconvex domain in $\C^n$ ($n > 1$) with real analytic simply connected boundary. Assume that $f$ is a nonconstant holomorphic mapping in a neighbourhood $U$ of a point $p \in \partial\Omega$ such that $f(U \cap \partial\Omega) \subset \partial\B^n$. Then $f$ extends to a biholomorphic mapping between $\Omega$ and $\B^n$.
\end{cor}

The proof uses the Reflection Principle. In our case the equation (\ref{Schwarz1}) has the form $(f,f) - 1 = 0$. Hence, one can  apply in (\ref{Schwarz3}) the Cramer rule (instead of the Implicit Function Theorem). It follows that $H$ in (\ref{Schwarz4}) is a rational function in $f$ and $X_jf$. This allows one to extend $f$ along the family $l_c$ of complex lines meromorphically but ``far enough". With this the proof can be completed as follows. First we complexify $\partial\Omega$ near a given point $p$ by a biholomorphic change of coordinates. Then we cut off a piece of $\partial\Omega$ by a real hyperplane parallel to the tangent plane at $p$. Next we extend $f$ meromorphically along a family of complex lines parallel to this hyperplane. Repeating this procedure we obtain a global meromorphic extension of $f$. The last step is to prove that this extension is, in fact, holomorphic.

A real analytic hypersurface $\Gamma$ is called {\it spherical} at a point $p\in \Gamma$ if in a neighbourhood of $p$ it is locally biholomorphic to an open piece of the real sphere $\partial\B^n$. It follows from Theorem~\ref{sphereTheo} that if a connected 
$\Gamma$ is spherical at one point, then it is spherical at every point.
Burns and Shnider~\cite{BuSh1} constructed the following  example. Let $\Gamma = \{z \in \C^2: y_2 = \vert z_1 \vert^2 \}$ (the unbounded sphere) and $\Gamma' = \{ z \in \C^2: \sin \ln \vert z_2 \vert^2 = 0, e^{-\pi} \le \vert z_2 \vert \le 1 \}$. Then the mapping $f(z) = (z_1/\sqrt{z_2}, \exp(i\ln z_2))$ with a suitable choice of a branch of $\ln z_2$, maps $\Gamma \setminus \{ 0 \}$ into $\Gamma'$ but does not extend even continuously to $z=0$.
In this example $\Gamma'$ is a compact real analytic spherical hypersurface which is not simply-connected. 

A result similar to Theorem~\ref{sphereTheo} holds if the sphere in the target space is replaced with an algebraic hypersurface,
see~\cite{Sha1}.

\begin{theorem}\label{t.ac}
Let $\Gamma$ be a connected essentially finite smooth real-analytic hypersurface in ${\mathbb C}^n$ and $p\in\G$. Let $\Gamma'$ be a compact strictly pseudoconvex real-algebraic hypersurface in ${\mathbb C}^n$. Let $f$ be a germ of a holomorphic mapping from $\Gamma$ to $\Gamma'$ defined at $p$. Then $f$ extends holomorphically along any path on $\Gamma$ with the extension
sending $\G$ to $\G'$.
\end{theorem}

In this result the hypersurface $\G$ is not assumed to be strictly pseudoconvex, although it follows from the proof that it is
pseudoconvex, and that the set of weakly pseudoconvex points of $\G$ consists precisely of the points where the extended map
degenerates. The proof is based on the technique of Segre varieties. The hypersurface $\Gamma$ is called essentially finite if the 
map $z \to Q_z$ is locally finite near every point of $\G$. The main idea is the holomorphic extension along Segre varieties:
from the properties of Segre varieties (see Section~\ref{s.geom}) one may conclude that the inclusion $f(Q_z) \subset Q'_{f(z)}$
must hold not only for points in the domain of $f$ but also for points $z$ with the property that $Q_z$ has a nonempty intersection
with the open set where $f$ is defined. This gives a holomorphic extension of $f$ to such points, which can be quite far away from the initial domain of $f$ no matter how small it is. An iterative procedure then gives extension of $f$ along any curve in $\Gamma$.
In particular, this technique can be used to give an alternative independent proof of Theorem~\ref{sphereTheo}.  In~\cite{KoSha},
 this technique was further refined to show that the germ of a biholomorphic mapping 
$f: \Gamma \to \Gamma'$ extends also across complex hypersurfaces that might be present in $\Gamma$. Here $\G'$ is either
a sphere or, more generally, any nondegenerate hyperquadric in $\mathbb C^n$.
 
Consider now the case of a real analytic hypersurface in the target domain. It turns out that unlike the spherical case, for 
nonspherical strictly pseudoconvex hypersurfaces the phenomenon of analytic continuation holds without any additional 
topological restrictions. This follows from  the following  result.

\begin{theorem}\label{t.cont2}
Let $\Gamma$, $\Gamma'$ be nonspherical real analytic strictly pseudoconvex hypersurfaces in $\C^n$, 
($n > 1$), and $U$ is a neighbourhood of a point $p \in \Gamma$, where the sets $\Gamma$, $\Gamma'$, $U$ and $\Gamma \cap U$ are connected and $\Gamma'$ is compact. Suppose that there exists a nonconstant holomorphic mapping $f:U \to \C^n$ such that 
$f(U \cap \Gamma) \subset \Gamma'$. Then $f$ continues analytically along any path in $\Gamma$ as a 
locally biholomorphic mapping.
\end{theorem}

In this form this theorem was proved in \cite{Pi55}. The proof is based on a careful analysis of the behaviour of Moser's chains on a nonspherical strictly pseudoconvex hypersurface. Vitushkin -Ezhov-Kruzhilin \cite{VEK} obtained a different proof of this result in a more general setting when $\Gamma$, $\Gamma'$ are nonspherical real analytic strictly pseudoconvex hypersurfaces in arbitrary 
$n$-dimensional complex manifolds, $n \ge 2$.
We refer the reader to the surveys \cite{Vi,Vi2} by Vitushkin for a comprehensive discussion and further results in this direction.

For nonalgebraic $\G'$, the problem of analytic continuation remains open in the presence of weakly pseudoconvex points in $\Gamma$, for example, 
it is not known if the map $f$ in Theorem~\ref{t.cont2} extends to a neighbourhood of an isolated weakly pseudoconvex point that $\Gamma$ might have. The difficulty is that the results of Moser's theory do not hold in general near points where the Levi form of $\G$ 
degenerates, and it is also not clear how to generalize the technique of analytic continuation along Segre varieties for nonalgebraic
target hypersurfaces. 

We conclude this section with a result by Nemirovski-Shafikov~\cite{NemSha1,NemSha2} on uniformization of strictly pseudoconvex domains.

\begin{theorem}\label{t.ns}
Let $\Omega$, $\Omega'$ be strictly pseudoconvex domains with real analytic boundaries. Then the universal coverings of $\Omega$ and $\Omega'$ are biholomorphically equivalent if and only if the boundaries of these domains are locally biholomorphically equivalent.
\end{theorem}

If the boundaries of $\Omega$ and $\Omega'$ are locally equivalent somewhere, then by Theorems~\ref{sphereTheo} and~\ref{t.cont2} the germ of the equivalence map extends as a locally biholomorphic map $f$ along any path in $\partial\Omega$, and hence along any path in a one-sided neighbourhood $V$ of $\partial\Omega$. This possibly multiple-valued map on $V$ 
extends to the envelope of holomorphy of $V$, which is~$\Omega$ by Hartogs' theorem. Kerner's 
theorem~\cite{Ker} states that the envelope of holomorphy of the universal covering $\widehat{V}$ of $V$ is the 
universal covering of the envelope of holomorphy of $V$. From this it follows that $f$ extends to a map 
$f: \widehat{\Omega} \to \Omega'$. In the nonspherical case, the final result can be deduced by repeating the argument 
for the inverse of the equivalence of $\partial \Omega$ and $\partial \Omega'$. In the spherical case an additional argument using invariant metrics in needed. The proof of Theorem~\ref{t.ns} in the other direction essentially follows the general scheme outlined in Section~\ref{s.refl}: 
the equivalence map $f: \hat\Omega \to \hat\Omega'$ is first extended smoothly to the boundary, and then
 the Reflections Principle (Theorem~\ref{SchwarzTheo1}) is applied.

\section{Weakly pseudoconvex domains}

\subsection{ Finite type and PSH peak functions.} Consider a smooth real hypersurface $\Gamma = \rho^{-1}(0)$ in $\C^n$. The following notion of the type of $\G$ at a point $p \in \Gamma$ is due to D'Angelo \cite{Dan1}. Denote by $O_p$ the space of germs of 
holomorphic mappings $h: (\C,0) \to \C^n$, $h(0) = p$. Denote by $\nu(h)$ the order of vanishing of $h - h(0)$ at the origin. Let 
also $\nu(\rho \circ h)$ denote the order of vanishing of the function $\rho \circ  h$ at the origin. Then the type of $\Gamma$ at $0$ is defined as
\begin{eqnarray}
\label{type1}
\tau(\Gamma,p) = \sup \{\nu(\rho \circ h)/\nu(h), h \in   O_p \} .
\end{eqnarray}
In general the function $p \mapsto \tau(\G,p)$ is not upper semicontinuous. Nevertheless, D'Angelo proved the following. 
Let $\Omega$ be a smoothly bounded pseudoconvex domain and let $q \in \partial\Omega$. Then there exists a neighbourhood 
$U$ of $q$ such that for each $p \in U \cap \partial\Omega$, 
\begin{eqnarray}
\label{type2}
\tau(\partial\Omega,p) \le \tau(\partial\Omega,q)^{n-1}/2^{n-2} .
\end{eqnarray}
A real analytic hypersurface is of finite type iff it contains no germs of complex analytic sets of positive dimension.
It is a result of Diederich-Fornaess \cite{DiFo2} that a compact real analytic subset of $\mathbb C^n$ contains no
nontrivial complex analytic subsets, and so every  bounded domain with real analytic boundary in $\C^n$ is of finite 
type at every boundary point.

The following characterization of finite type is due to Fornaess-Sibony \cite{FoSi, Si3}:
\begin{prop}
\label{type3}
Let $\Omega$ be a smoothly bounded domain in $\C^n$, $p \in \partial\Omega$. Assume there exist a function $\phi_p \in C^0(\overline\Omega)$ plurisubharmonic in $\Omega$ and constants $C > 0$, $\lambda > 0$, $k > 0$ such that 
\begin{eqnarray}
\label{type4}
-C\vert z - p \vert^{\lambda} \le \phi(z) \le - \vert z \vert^{2k\lambda}, \ \ z \in \overline\Omega .
\end{eqnarray}
Then $\partial\Omega$ is of type less than $2k$ at $p$.
\end{prop}
The function $\phi_p$ above is called a plurisubharmonic {\it barrier} (at $p$). If $k=1$, the existence of a barrier is equivalent 
to strict pseudoconvexity at $p$ (see \cite{Si2}). 

We say that a boundary point $q \in \partial\Omega$ satisfies {\it the barrier property} if  there exists a neighbourhood $U$ of 
$q$ such that every point  $p \in U \cap \partial\Omega$ admits a barrier function (with $k$ and $\lambda$ independent of~$p$).

\begin{theorem}
\label{type5}
Let $\Omega$ be a domain in $\C^n$ with a smooth pseudoconvex boundary of finite type in a neighbourhood $V$ of a point $q \in \partial\Omega$.
Then $q$ satisfies the barrier property.
\end{theorem}

When $n=2$ or when $V \cap \partial\Omega$ is convex, this result is due to Fornaess-Sibony \cite{FoSi}. The real analytic 
case is due to Diederich-Fornaess \cite{DiFo4}. Finally, the general case was treated by Cho \cite{Cho}.

Diederich-Fornaess \cite{DiFo3} proposed the use of barrier functions in order to obtain lower bounds for 
the Kobayashi-Royden metric. This approach works for a wide class of domains. As an application we present the following 
result obtained in \cite{Su2}.

\begin{prop}
\label{ExtTheo2}
Let $\Omega$ and $\Omega'$ be domains  in $\C^n$ whose boundaries are $C^2$-smooth near some points 
$p \in \partial\Omega$ and $p' \in \partial\Omega'$ and satisfy the barrier property at these points. Let  
$f:\Omega \to \Omega'$ be  a holomorphic mapping such that for some neighbourhood $U$ of $p$ the cluster set 
$C_\Omega(f; U \cap \partial\Omega)$ doest not intersect $\Omega'$. Assume also that the cluster set $C_\Omega(f;p)$ 
contains the point $p'$. Then $f$ extends to a neighbourhood of $p$ in $\partial\Omega$ as a H\"older-continuous 
mapping.
\end{prop}

Various results concerning continuous extension of holomorphic mappings are also obtained in \cite{Ber1,Ber2,DiFo3,FoLow,FoRo}.

\subsection{$\overline{\partial}$-approach to boundary regularity} 
Let $\Omega$ be a bounded domain in $\C^n$. Consider the hermitian Hilbert space $L^2(\Omega)$ equipped with the standard hermitian product $(\bullet,\bullet)_{L^2(\Omega)}$. Then ${\mathcal O}(\Omega) \cap L^2(\Omega)$ is a closed subspace in $L^2(\Omega)$ and hence is itself a Hilbert space. Fix a point $p \in \Omega$. The evaluation map 
$$l_p:{\mathcal O}(\Omega) \cap L^2(\Omega) \to \C, \,\,\, l_p: h \mapsto h(p),$$
is a bounded linear functional on ${\mathcal O}(\Omega) \cap L^2(\Omega)$. By the Riesz Representation Theorem there exists a unique element in ${\mathcal O}(\Omega) \cap L^2(\Omega)$, which is denoted by $K_\Omega(\bullet,p)$, such that
$$h(p) = l_p(h) = (h,K_\Omega(\bullet,p))_{L^2(\Omega)}$$
for all $h \in {\mathcal O}(\Omega) \cap L^2(\Omega)$. The function $K_\Omega: \Omega \times\Omega \to \C$ is called the Bergman kernel for $\Omega$.
By this definition the function $z \mapsto K_\Omega(z,p)$ is in $L^2(\Omega)$ for evey $p \in \Omega$. Furthermore, $K_\Omega(z,p) = \overline{K_\Omega(p,z)}$ and the function $(z,w) \mapsto K_\Omega(z,\overline{w})$ is holomorphic on $\Omega \times \Omega$. The orthogonal projection operator 
$$P_\Omega: L^2(\Omega) \to {\mathcal O}(\Omega) \cap L^2(\Omega)$$
is called the Bergman projection. One has $P_\Omega(h) = (h,K_\Omega)_{L^2(\Omega)}$. The following transformation rules (see, for example, \cite{Ra}) play a key role in application of the Bergman kernel and the Bergman projection to holomorphic mappings.
\begin{theorem}
\label{Bergman1}
Let $f:\Omega_1 \to \Omega_2$ be a biholomorphic mapping between bounded domains in $\C^n$. Then
\begin{eqnarray*}
& &K_{\Omega_1}(p,z) = J_f(p) K_{\Omega_2}(f(p),f(z))\overline{J_f(z)},\\
& &P_{\Omega_1}(J_f h \circ f) = J_f (P_{\Omega_2}(h) \circ f),
\end{eqnarray*}
for all $h \in L^2(\Omega_2)$. Here $J_f$ denotes the determinant of Jacobian matrix of $f$.
\end{theorem}
The above transformation rule for the Bergman projection (but not for the kernel) remains true also for proper holomorphic mappings.

A smoothly bounded pseudoconvex domain $\Omega$ is said to satisfy {\it Condition R} if 
$$P_\Omega(C^\infty(\overline{\Omega})) \subset C^\infty(\overline{\Omega}).$$
The following result is due to Bell-Catlin \cite{BeCa1} and Diederich-Fornaess \cite{DiFo}.

\begin{theorem}\label{Bergman2}
Let $f:\Omega_1 \to \Omega_2$, be a proper holomorphic mapping between smoothly bounded pseudoconvex domains. Suppose that $\Omega_1$ satisfies Condition R. Then $f$ extends as a $C^\infty$-smooth mapping on $\overline{\Omega_1}$.
\end{theorem}

A version of this theorem for $CR$ mappings between the boundaries of domains was obtained by Bell-Catlin~\cite{BeCa2}.

A general approach to verify Condition R for a prescribed class of domains relies on the $\overline\partial$-Neumann problem. Let $\Omega$ be a smoothly bounded pseudoconvex domain in $\C^n$. Let $g$ be a $\overline\partial$-closed (0,1)-differential form with coefficients of class $L^2(\Omega)$, i.e., $g \in L^2_{(0,1)}(\Omega)$. The $\overline\partial$-Neumann problem consists of determining the regularity of the solution $u$ to the equation $\overline\partial u = g$ which is orthogonal to the kernel of the operator $\overline\partial$, that is, to the class ${\mathcal O}(\Omega) \cap L^2(\Omega)$. This solution is called the canonical solution. The operator 
$$N_\Omega: L^2_{(0,1)}(\Omega) \to L^2(\Omega), \,\,\, N_\Omega: g \mapsto u,$$
is called the $\overline\partial$-Neumann operator on $\Omega$. The relation between the Bergman projection and the $\overline\partial$-Neumann operator is given by Kohn's formula \cite{Ko}:
$$P_\Omega = Id - \overline\partial^* N_\Omega \overline\partial.$$
Thus, if the $\overline\partial$-Neumann operator is globally regular, i.e., maps the space $C^\infty(\overline\Omega)$  to itself, 
then Condition R holds. Regularity of the $\overline\partial$-Neumann operator has been an active area of research and lead to
the development of many important technical tools, see Catlin~\cite{Ca1, Ca2, Ca3}. In particular, it is known that the existence of plurisubharmonic barriers (\ref{type4}) implies the regularity of the 
$\overline\partial$-problem (see \cite{Ca4, Si3}).

A smoothly bounded domain $\Omega \subset \C^n$ admits a defining function which is plurisubharmonic along the boundary 
if there exists a smooth defining function of $\Omega$ whose Levi form is positive semi-definite for all vectors at each boundary point (this condition is stronger than the pseudoconvexity which requires that the Levi form is positive semi-definite on the holomorphic tangent space). Boas and Straube \cite{Boas1} proved that if $\Omega$ admits a defining function which is plurisubharmonic 
along the boundary, then it satisfies Condition R. In particular, every smoothly bounded convex domain satisfies Condition R.

Finally, it was shown by Christ~\cite{Chr} that for the worm-domain of Diederich-Fornaess~\cite{DF0}, which is smooth and pseudoconvex, Condition~R does not hold.
This shows limitations of this approach and raises an important question of finding sufficient and necessary
conditions for the regularity of the Bergman projection.

\section{Automorphism groups and proper self-mappings}

The geometry of the boundary of a domain $\Omega$ in $\mathbb C^n$ influences the structure of the group of biholomorphic automorphisms of $\Omega$. In its turn, the automorphism group ${\rm Aut}(\Omega)$ may completely characterize the 
domain $\Omega$. In this section we discuss some results in this direction.

\subsection{Automorphisms of strictly pseudoconvex domains}
We begin with the following result, generally known in the literature as the Wong-Rosay theorem.

\begin{theorem}
\label{WRTheo1}
Let $\Omega$ be a strictly pseudoconvex domain in $\C^n$. Assume that ${\rm Aut}(\Omega)$ is not compact.
Then $\Omega$ is biholomorphic to the unit ball $\B$.
\end{theorem}
Since ${\rm Aut}(\Omega)$ is not compact, there exists a sequence $\{f^k\}$ in ${\rm Aut}(\Omega)$ which converges 
uniformly on every compact subset of $\Omega$ to a boundary point $q \in \partial\Omega$.
Fix a point $p \in \Omega$ and set $q^k = f^k(p)$. Let $D^k$ be the scaling sequence for $\{q^k\}$. Then 
the sequence $F^k = D^k \circ f^k$ converges to a biholomorphic mapping from $\Omega$ to $\B$.

Theorem \ref{WRTheo1} was established by Webster \cite{We2} under the additional assumption that the group 
${\rm Aut}(\Omega)$ has positive dimension. In full generality this result was obtained by Burns-Shnider~\cite{BuSh2} 
using the Chern-Moser theory. A more elementary approach based on invariant metrics is due to Wong~\cite{Wo} and 
Rosay~\cite{Ro}. The short proof presented above was given in~\cite{Pi7,Pi8}.

There are several other local versions of this result, for example \cite{Ef}: 

\begin{theorem}
\label{WRTheo2}
Let $\Omega$ be a domain (not necessarily bounded) in $\C^n$ and let $\partial\Omega$ be $C^2$-smooth strictly 
pseudoconvex in a neighbourhood of a point $q \in \partial\Omega$. Suppose that there exist a sequence $\{f^k\}$ in 
${\rm Aut}(\Omega)$ and a point $p \in \Omega$ such that $f^k(p) \to q$ as $k \to \infty$. Then $\Omega$ is 
biholomorphic to the unit ball $\B^n$.
\end{theorem}

It is also interesting to consider the inverse question: which groups can be realized as an automorphism group of a domain in 
$\C^n$? The following results are due to Winkelman \cite{Wi1,Wi2}:
\begin{theorem}
\label{WRTheo3}
Let $G$ be a (finite or infinite) countable group. Then there exists a (connected) Riemann surface $M$ 
such that $G$ is isomorphic to ${\rm Aut}(M)$.
\end{theorem}

\begin{theorem}
\label{WrTheo3}
Let $G$ be a connected (real) Lie group. Then there exists a Stein, complete hyperbolic complex manifold $M$ on which $G$ acts effectively, freely, properly and with totally real orbits such that $G$ is isomorphic to ${\rm Aut}(M)$.
\end{theorem}

Note that any compact real Lie group can be realized as an automorphism group of a strictly pseudoconvex 
domain~\cite{BeDa,SaZ}.

\subsection{Domains with large automorphism groups}
Investigation of weakly pseudoconvex domains with large automorphism groups was initiated by Greene-Krantz~\cite{GK}.
The following result is due to Bedford-Pinchuk \cite{BePi1, BePi4}.

\begin{theorem}\label{t.b-p}
Let $\Omega$ be a bounded pseudoconvex domain with real analytic boundary in $\C^2$. Suppose that ${\rm Aut}(\Omega)$ is not compact. Then $\Omega$ is biholomorphic to a domain of the form 
\begin{eqnarray}
\label{BePi1}
\{ z \in \C^2: \vert z_1 \vert^2 + \vert z_2 \vert^{2m} < 1 \}
\end{eqnarray}
for some positive integer $m$.
\end{theorem}

The same result also holds if $\Omega$ is a smoothly bounded pseudoconvex domain of finite type. Furthermore, the assumption of 
pseudoconvexity (if the boundary is real analytic) can be dropped, see \cite{BePi2}.

We outline the proof of Theorem~\ref{t.b-p}. Since ${\rm Aut}(\Omega)$ is noncompact, there exists a point $a \in \Omega$ 
and a sequence of 
automorphisms $\{f^j\}$ such that the sequence  $q^k = f^k(a)$ converges to a boundary point $q \in \partial\Omega$. 
Then the sequence $(f^k)$  converges uniformly on compact subsets of $\Omega$ to a constant map $f^0  \equiv q$. 
Applying the scaling along the sequence $(q^k)$ one can prove that $\Omega$ is equivalent to a domain of the form 
$D = \{ 2x_2 + P(z_1,\overline{z}_1) = 0 \}$, where $P$ is a nonzero real polynomial. 
Note that the proof is more delicate than in the strictly pseudoconvex case and is based on precise estimates of the 
Kobayashi-Royden  metric in pseudoconvex domains of finite type in $\C^2$. These estimates were obtained by 
Catlin~\cite{Ca5}; a geometric proof of his result based on the scaling method was given by Berteloot~\cite{Ber5}. 
The one-parameter group $L^t(z_1,z_2) = (z_1,z_2 + it)$ acts on the domain~$D$. The biholomorphism $f:D \to \Omega$ 
defines a real one-parameter group of automorphisms $h^t = f \circ L^t \circ f^{-1}$. 

The second step is to prove that the group $(h^t)$ is parabolic, that is, there exists a point $p\in \Omega$ (called a parabolic 
point) such that 
$$\lim_{t\to-\infty}h^t(z) = \lim_{t \to \infty}h^t(z) = p.$$
The proof also uses the estimates of the Kobayashi metric.

The next step is to study the holomorphic vector field $X = (X_1,X_2)$ generating the parabolic subgroup $(h^t)$. 
This vector field is tangent to $\partial\Omega$, that is, 
$$\Re \left (\frac{\partial \rho}{\partial z_1}X_1 + \frac{\partial \rho}{\partial z_2}X_2\right) = 0 .$$
This condition leads to a rather precise description of the jet of $\partial\Omega$ at 
a parabolic point $p$, and this can be used to conclude the proof.

In the local case  Verma \cite{Ve} obtained the following classification result.

\begin{theorem}
Let $\Omega$ be a bounded domain in $\C^2$. Suppose that there exists a point $p \in \Omega$ and a sequence $\{ \phi_j\} \in {\rm Aut}(\Omega)$ such that $\{ \phi_j(p)\}$ converges to $p_\infty \in \partial\Omega$. Assume that the boundary of $\Omega$ is real analytic and of finite type near $p_\infty$. Then exactly one of the following cases holds:
\begin{itemize}
\item[(i)] If $\dim {\rm Aut}(\Omega) = 2$ then either
\begin{itemize}
\item[(a)] $\Omega$ is biholomorphic to $\Omega_1 = \{ z \in \C^2: 2 \Re z_2 + P_1(\Re z_1) < 0 \}$ where $P_1(\Re z_1)$ is a polynomial that depends on $\Re z_1$, or
\item[(b)] $\Omega$ is biholomorphic to $\Omega_2 = \{ z \in \C^2: 2\Re z_2 + P_2(\vert z_1 \vert^2) < 0 \}$ where $P_2(\vert z_1 \vert^2)$ is a homogeneous polynomial that depends on $\vert z_1 \vert^2$, or
\item[(c)] $\Omega$ is biholomorphic to $\Omega_3 = \{ z \in \C^2: 2\Re z_2 + P_{2m}(z_1,\overline{z}_1) < 0 \}$ where $P_{2m}(z_1,\overline{z}_1)$ is a homogeneous polynomial of degree $2m$ without harmonic terms.
\end{itemize}
\item[(ii)] If $\dim {\rm Aut}(\Omega) = 3$ then $\Omega$ is biholomorphic to $\Omega_4 = \{ z \in \C^2: 2\Re z_2 + (\Re z_1)^{2m} < 0 \}$ for some integer $m \ge 2$.
\item[(iii)] If $\dim {\rm Aut} (\Omega) = 4$ then $\Omega$ is biholomorphic to $\Omega_5 = \{ z \in \C^2: \vert z_1 \vert^2  + \vert z_1 \vert^{2m} < 0 \}$ for some integer $m \ge 2$. 
\item[(iv)] If $\dim {\rm Aut}(\Omega) = 8$ then $\Omega$ is biholomorphic to $\B^2$.
\end{itemize}
The dimensions $0,1,5,6,7$ cannot occur with $\Omega$ as above.
\end{theorem}

In higher dimensions the situation is more complicated. We assign to the variables $z_1,...,z_n$ the weights 
$\delta_1,...,\delta_n$, where $\delta_j = (2m_j)^{-1}$ for $m_j$  a positive integer.
If $J = (j_1,...,j_n)$ and $K=(k_1,...,k_n)$ are the multi-indices, we set 
$wt(J) = j_1\delta_1 + ...+j_n\delta_n$ and $wt(z^J\overline{z}^K) = wt(J) + wt(K)$. 
We consider real polynomials of the form 
\begin{eqnarray}
\label{BePi2}
p(z,\overline{z}) = \sum_{wt J = wt K = 1/2} a_{JK}z^J\overline{z}^K .
\end{eqnarray}
The reality of $p$ is equivalent to $a_{JK} = \overline{{a}_{KJ}}$. The balance of the weights $wt(J) = wt(K)$ implies 
that the domain
\begin{eqnarray}
\label{BePi3}
G = \{ (w,z_1,...,z_n) \in \C \times \C^n: \vert w \vert^2 + p(z,\overline{z}) < 1 \}
\end{eqnarray}
is invariant under the action of the real torus
\begin{eqnarray}\label{BePi4.1}
(\phi,\theta) \mapsto (e^{i\phi}w,e^{i\delta_1\theta}z_1,...,e^{i\delta_n\theta}z_n) .
\end{eqnarray}
The weighted homogeneity of $p$ implies that the Cayley-type transform $(w,z) \mapsto (w^*,z^*)$, defined by 
\begin{eqnarray}\label{BePi4.2}
w = (1 - iw^*/4)(1 + iw^*/4)^{-1},\,\,\,z_j = z_j^*(1 + iw^*/4)^{-2\delta_j} ,
\end{eqnarray}
maps $G$ biholomorphically onto the domain 
\begin{eqnarray}
\label{BePi5}
D = \{ (w,z_1,...z_n): \C \times \C^n: \Im w + p(z,\overline{z}) < 0 \}.
\end{eqnarray}
The latter is an unbounded realization of $G$. Note that $D$ is invariant under the translation along the
$\Re w$-direction. Since $p$ is homogeneous, the domain $D$ is invariant with respect 
to the family of anisotropic dilations. Hence the dimension of ${\rm Aut}(D)$ is at least $4$.

\begin{theorem}
\label{BePi6}
Let $\Omega \in \C^{n+1}$ be a smoothly bounded convex domain of finite type. If ${\rm Aut}(\Omega)$ is 
noncompact then $\Omega$ is equivalent to the domain of the form (\ref{BePi3}).
\end{theorem}
This result is obtained in \cite{BePi3}. The scaling method in a convex domain $\Omega$ (not necessarily of finite type) relies 
on the estimates of the Kobayashi-Royden metric, which also have 
other applications. Denote by $L(a,V)$ the complex line passing through a point $a \in \Omega$ in the direction of a  vector $V$.  Define $\delta(a,V)$ to be the Euclidean distance from $a$ to $L(a,V) \cap \partial\Omega$. Then the following estimate holds \cite{BePi3,Fra}:
\begin{eqnarray}
\label{BePiEst}
\frac{\vert V \vert }{ 2 \delta(a,V)} \le F_\Omega(a,V) \le \frac{\vert V \vert }{ \delta(a,V)} .
\end{eqnarray}
Using this estimate and his version of the scaling method Frankel \cite{Fra} proved the following

\begin{theorem}
Suppose that $\Omega \subset \mathbb C^n$ is a bounded convex domain and that there exists a discrete subgroup of ${\rm Aut}(\Omega)$ 
which acts properly discontinuously, freely, and cocompactly on $\Omega$. Then $\Omega$ is a bounded symmetric domain.
\end{theorem}

Recently, Zimmer \cite{Zi} proposed a new approach to the problem of classification of convex domains with large automorphism groups. Let $\Omega$ be a domain in $\C^n$. {\it The limit set} of $\Omega$ is the set of points $z \in \partial\Omega$  
for which there exist some $p \in \Omega$ and some sequence $\phi_k \in {\rm Aut}(\Omega)$ such that $\phi_k(p) \to z$. 
If ${\rm Aut}(\Omega)$ is noncompact, the limit set is not empty. If $\Omega$ is a bounded convex domain with $C^1$-smooth boundary, {\it the closed complex face} of a point $z \in \partial\Omega$ is the closed set $\partial\Omega \cap H_z(\partial\Omega)$. The main result of \cite{Zi} is the following

\begin{theorem}
Suppose $\Omega$ is a bounded convex domain with $C^\infty$-smooth boundary. Then the following are equivalent:
\begin{itemize}
\item[(1)] The limit set of $\Omega$ intersects at least two closed complex faces of $\partial\Omega$.
\item[(2)] $\Omega$ is biholomorphic to (\ref{BePi3}).
\end{itemize}
\end{theorem}

Notice that there is no finite type assumption in this theorem. The main new tool used by Zimmer is the theory of Gromov 
hyperbolic metric spaces.

Suppose that $(X,d)$ is a metric space. Let $I \subset \R$ be an interval. A curve $\sigma:I \to X$ is called {\it geodesic} if 
$d(\sigma(t_1),\sigma(t_2)) = \vert t_1 - t_2 \vert$ for all $t_1,t_2 \in I$. A {\it geodesic triangle} is a choice of 3 points in $X$ and geodesic segments connecting these points. A geodesic triangle is said to be $\delta$-{\it thin} if any point on any of the sides of the triangle is within distance $\delta$ of the other two sides. A proper geodesic metric space $(X,d)$ is called $\delta$-{\it hyperbolic} if every geodesic triangle is $\delta$-thin. If $(X,d)$ is $\delta$-hyperbolic for some $\delta \ge 0$, then $(X,d)$ is called {\it Gromov hyperbolic}.   Zimmer's approach uses the following result  established  in \cite{Zi1}:
\begin{theorem}
\label{GrKoTh}
Suppose $\Omega \subset \C^n$ is a bounded convex domain with smooth boundary. Then the following are equivalent:
\begin{itemize}
\item[(1)] $\Omega$ has a finite type in the sense of D'Angelo,
\item[(2)] $(\Omega,d_\Omega)$ is Gromov hyperbolic, where $d_\Omega$ is the Kobayashi distance on $\Omega$.
\end{itemize}
\end{theorem}
Gromov hyperbolicity of strictly pseudoconvex domains was established by Ballog-Bonk \cite{BaBo} using estimates of the 
Kobayashi-Royden metric (see Proposition \ref{KobProp1}). Further results in this direction are obtained recently by Bracci-Gaussier \cite{BrGa}.

The condition of convexity is crucially used in the proofs of the above results. It is not known whether any bounded pseudoconvex 
domain with smooth boundary of finite type and noncompact automorphisms group in $\C^n$, $n > 2$, is equivalent to 
(\ref{BePi3}). The question remains open even for domains with real algebraic boundaries.

\subsection{Proper self-mappings}
Finally, we discuss some progress in the direction originated from Alexander's Theorem \ref{AlRu}: a proper holomorphic self-map of $\B^n$, $n > 1$, is a biholomorphism. Its generalization to the case of strictly pseudoconvex domain (see Corollary \ref{AlPinCor}) is based on Theorem~\ref{AlPinTheo}. However, the condition of strict pseudoconvexity is crucial for Theorem \ref{AlPinTheo}. Indeed, the map $f:(z_1,z_2) \to (z_1^2,z_2)$ takes the domain $\{ z\in \C^2: \vert z_1\vert^4 + \vert z_2 \vert^2\}$ properly on $\B^2$ but its  critical locus is not empty. This is a serious obstacle for applications of the scaling method and for this reason the analogs of Alexander's theorem are currently established only for special classes of domains.

One of the most general results in this direction was obtained by Bedford \cite{Bed1}.
\begin{theorem}

\label{BedTheo}
Let $\Omega$ be a bounded pseudoconvex domain with real analytic boundary in $\C^n$, $n \ge 2$. Then every proper holomorphic self-mapping
$f:\Omega \to \Omega$ is a biholomorphism.
\end{theorem}

The proof is based on a careful analysis of the branch locus of a proper holomorphic mapping from a pseudoconvex domain with 
real analytic boundary. Notice that the assumption of real analyticity is crucially used here. To the best of our knowledge it is not 
known if the analog of Theorem~\ref{BedTheo} remains true
for pseudoconvex domains with smooth boundary of finite type in the sense of D'Angelo. Certain results of this type are obtained for domains which admit some symmetries.

A domain $\Omega$ is said to be quasi-regular if there exist integers $p$ and $q$, $p + q \ge 1$, such that whenever $(z,w) \in \Omega$, $(e^{ip\theta},e^{iq\theta}) \in \Omega$ for $\theta \in [0,2\pi]$. Thus, if $p=q=1$, the domain $\Omega$ is circular; when $p=0$ or $q=0$, $\Omega$ is a Hartogs domain. The following result is obtained in \cite{CoPaSu1,CoPaSu2}:
\begin{theorem}
Let $\Omega$ be a smoothly bounded pseudoconvex quasi-circular domain of finite type in $\C^2$. Then every proper holomorphic self-map of $\Omega$ is a biholomorphism.
\end{theorem}
The proof uses the scaling method (for the study of the branch locus of the map) and arguments from holomorphic dynamics. 

Assumptions on regularity of the boundary can be weakened for domains with additional symmetries. For every $a \in \C^n$ denote by $L_a:\C^n \to \C^n$
the linear map $T_az = (a_1z_1,...,a_nz_n)$. Recall that a domain $\Omega \subset \C^n$ is called a Reinhardt domain (resp. complete) if $T_a(\Omega) = \Omega$ for every $a$ with $\vert a_j \vert = 1$ (resp. $\vert a_j \vert \le 1$), $1 \le j \le n$. 
We note that a complete description of automorphisms of a wide class of hyperbolic Reinhardt domains was obtained by 
Kruzhilin~\cite{Kru1}.

The following result was established by Berteloot-Pinchuk \cite{BerPi}:

\begin{theorem}
\label{BePiTheo}
Among bounded, complete, Reinhardt domains in $\C^2$, the bidiscs are the only ones that admit proper holomorphic 
self-mappings that are not automorphisms.
\end{theorem}

This work also contains a detailed description of proper holomorphic maps between complete Reinhardt domains. The general case of Reinhardt domains in $\C^2$ (not necessarily complete) was considered by Isaev-Kruzhilin \cite{IsKru}. They obtained a complete description of proper holomorphic mappings and classified all Reinhardt domains in $\C^2$ admitting proper holomorphic self-maps which are not biholomorphisms. A partial generalization of Theorem \ref{BePiTheo} to higher dimensions is obtained by 
Berteloot~\cite{Ber4}.

Proper holomorphic mappings between the classical Cartan domains and a wide class of Siegel domains were studied by Tumanov-Henkin \cite{TuHe1,TuHe2} 
and Henkin-Novikov \cite{HeNo}. We present here one of their results:
\begin{theorem}
Let $\Omega \subset \C^n$, $n > 1$, be an irreducible bounded symmetric domain. Then every proper holomorphic self-map $f:\Omega \to \Omega$ is an automorphism of $\Omega$.
\end{theorem}

\section{Proper holomorphic mappings between real analytic domains} 
The goal of this section is to present the following results obtained by Diederich-Pinchuk \cite{DiPi1,DiPi2003}.

\begin{theorem}\label{DP1}
Let $f: \Omega \to \Omega'$ be a proper holomorphic mapping between two bounded domains in $\C^n$ with real analytic 
boundaries. Suppose that at least one of the following conditions holds:
\begin{itemize}
\item[(a)] $n=2$;
\item[(b)] $f$ extends continuously on $\overline{\Omega}$.
\end{itemize}
Then $f$ extends holomorphically to a neighbourhood of $\overline{\Omega}$.
\end{theorem}

When the map $f$ is assumed to be a biholomorphism and to extend smoothly to the boundary of $\Omega$, then this result was obtained by 
Baouendi-Jacobowitz-Tr\`eves~\cite{BaJaTr}. For pseudoconvex domains this was proved in any dimension and 
without the assumption of boundary continuity 
by Diederich-Fornaess \cite{DiFo4} and Baouendi-Rothschild \cite{BaRo}. In that case pseudoconvex boundaries are 
automatically of finite type and Condition R holds. Therefore, a proper holomorphic map $f$ extends smoothly on 
$\overline{\Omega}$ by Theorem~\ref{Bergman2}. Case (a) was also considered by Huang \cite{Huang} under an 
additional assumption that $f$ is continuous on $\overline{\Omega}$.

Part (b) follows from a more general result.

\begin{theorem}
\label{DP2}
Let $\Gamma \subset \Omega$ (resp. $\Gamma' \subset \Omega'$) be a real analytic closed hypersurface of finite type  and let $f:\Gamma \to \Gamma'$ be a continuous CR mapping. Then $f$ extends holomorphically to a neighbourhood of $\Gamma$.
\end{theorem}

The proof of these results consists of two major parts:

(1) One proves that $f$ extends as a proper holomorphic correspondence to a neighbourhood of~$\partial\Omega$.

(2) One proves that if $f$ extends as a proper holomorphic correspondence then it extends as a holomorphic mapping to a neighbourhood of~$\partial\Omega$.

\medskip

The main tool in the proof is the invariance property of Segre varieties associated with the real analytic hypersurfaces. Their 
behaviour near Levi degenerate points of the boundary requires a more subtle analysis. We describe now the main steps.

Let $\Omega$ be a bounded domain in $\C^n$ with real analytic boundary. There exist a neighbourhood $W$ of $\partial\Omega$ and a real analytic function $\rho:W \to \R$ such that $D \cap W = \{ z \in W: \rho(z) < 0 \}$ and $d\rho(z) \neq 0$ for all $z \in \partial\Omega$ (a global defining function). Its complexification $\rho(z,\overline{w})$ is defined on a suitable neighbourhood $V \subset \C^{2n}$ of the diagonal $\Delta \subset W \times W$ and is holomorphic in $z$ and antiholomorphic in $w$. For points $z \in \C^{n}$ we use the notation $z = ('z,z_n) \in \C^{n-1} \times \C$.

Let $z^0 \in \partial\Omega$. A local holomorphic coordinate system centred at $z^0$ is called standard if the defining function $\rho$ can be written in these coordinates in the form $\rho(z) = 2 x_n + o(\vert z \vert)$. A pair of open neighbourhoods $U_1 \subset U_2$ (with $\overline{U_1} \subset U_2$) is called a standard pair of neighbourhoods of $z^0$ if it has the following properties:

\begin{itemize}
\item[(a)] With respect to a suitable standard coordinate system at $z^0$ one has $U_2 = {'U_2} \times U_{2n}$ with $'U_2$ 
being an open neighbourhood of $0 \in \C^{n-1}$,  and $U_{2n}$ an open neighbourhood on the $z_n$-axis.
\item[(b)] The complexification $\rho(z,\overline{w})$ is well defined on $U_2 \times U_1$ so that for each $w \in U_1$ the Segre variety $Q_w = \{ z \in U_2: \rho(z,\overline{w}) = 0 \}$ is well-defined.
\item[(c)] $Q_w$ can be written as a graph. This means that there exists a holomorphic function $h_w('z)$ on $'U_2$ (depending antiholomorphically on $w$), such that
\begin{eqnarray}
\label{e.DP1}
Q_w = \{ ('z,z_n) \in U_2: z_2 = h_w('z) \} .
\end{eqnarray}
\end{itemize}
 
 Note that every point $z^0 \in \partial\Omega$ admits a family of standard pairs of neighbourhoods such that the 
 corresponding  $U_2$ form a neighbourhood basis of $z^0$. With this notation the function $h('z,\overline{w}):= h_w('z)$ 
 can be written as a power series $h('z,\overline{w}) = \sum_j \lambda_j(\overline{w})z^j$ with  coefficients $\lambda_j$ antiholomorphic on $U_1$. There exists an integer $N$ (depending only on $\partial\Omega$) such that for all $z^0 \in \partial\Omega$ and any standard pair of neighbourhoods $U_1 \subset U_2$ of $z^0$, the coefficients $\{ \lambda_j: \vert j \vert \le N \}$ uniquely determine $Q_w$. This allows us to define the structure of a finite dimensional complex variety on 
 the family of all Segre varieties so that the maps
 \begin{eqnarray}\label{e.DP2}
 \lambda: U_1 \ni w \mapsto Q_w
 \end{eqnarray}
 are finite antiholomorphic branched coverings. For any point $w \in W$, with $W$ being a sufficiently small open neighbourhood of 
 $\partial\Omega$, one has the following: the complex line $l_w$ through $w$ containing the real line passing through $w$ and orthogonal to $\partial\Omega$ intersects the Segre variety $Q_w$ at exactly one point ${}^sw$, called the symmetric point of 
 $w$. For $w \in W \setminus \overline\Omega$ one always has ${}^sw \in \Omega$. The connected component of $Q_w \cap \Omega$ containing ${}^sw$ is denoted by ${}^sQ_w$ and is called the symmetric component.

The second important technical tool is provided by holomorphic correspondences. Let $U$, $U'$ be open subsets of $\C^n$. A proper holomorphic correspondence is a closed complex analytic subset $F \subset U \times U'$ of pure dimension $n$ such that the canonical projection $\pi:F \to U$ is proper. The correspondence $F$ is called irreducible if $F \subset U \times U'$ is irreducible as an analytic set (see~\cite{Ch1} for generalities on complex analytic sets). Let $\Omega$, $\Omega'$ be bounded domains in $\C^n$ and $z^0 \in \partial\Omega$ be a boundary point. We say that $f$ extends as a proper holomorphic correspondence to a neighbourhood $U$ of $z^0$ if there exist an open set $U' \subset \C^n$ and an irreducible proper holomorphic correspondence $F \subset U \times U'$ such that 
$$\Gamma_f \cap\{ (\Omega \cap U) \times \Omega' \} \subset F,$$
where $\Gamma_f$ denotes the graph of $f$. 

One can view this as an extension of $f$ as a multiple-valued map. Indeed, a correspondence $F$ assigns to each point $z \in U$ a 
finite number of points in the target space, namely, the set ${\hat F(z)} := \pi'(\pi^{-1}(z))$, where $\pi'$ denotes the projection 
of $F$ to $U'$. 

Let $f:\Omega \to \Omega'$ be a proper holomorphic mapping between two bounded domains with real analytic boundaries in $\C^n$. Suppose that $f$ extends as a correspondence $F$ to a neighbourhood of the point $z^0 \in \partial\Omega$. Choose standard coordinates such that $z^0 = 0$, $f(z^0) = 0$ and the standard neighbourhoods $U_j$,(resp. $U_j'$), $j=1,2$. Then we have the following invariance property for the Segre varieties under $\hat F$:
\begin{prop}
\label{DP3}
For every $(w,w') \in F \cap (U_1 \times U_1')$, the inclusion $\hat{F}(Q_w) \subset Q_{w'}'$ holds.
\end{prop}

Now we can explain how to construct a holomorphic correspondence which extends the graph of $f$. 
For $\zeta \in Q_w$ we denote by ${}_\zeta Q_w$ the germ of $Q_w$ at $\zeta$.
For every point $z^0 \in \partial\Omega$ in a standard coordinate system, a standard pair of neighbourhoods 
$U_1 \subset U_2$, and a suitably chosen open neighbourhood $U'$ of $\partial\Omega'$ we define
\begin{eqnarray}
\label{DP4}
V:= \{ (w,w') \in (U_1 \setminus \overline\Omega) \times (U' \setminus \overline{\Omega'}): {}_{{}^sw'}Q_{w'}' \subset f(Q_w \cap \Omega) \} .
\end{eqnarray}
The important step is to show that $V$ extends as an $n$-dimensional analytic set to a full neighbourhood of $(0,0)$,
which, in fact, is the extension of the graph of $f$. This is not obvious because it requires the properness of the 
projection of $V$ to  $U_1 \setminus \Omega$.

The second part of the proof is given by the following 

\begin{theorem}
Let $\Omega$, $\Omega' \subset \C^n$ be bounded domains with real analytic boundaries and $f:\Omega \to \Omega'$ be a proper holomorphic mapping that extends as a holomorphic correspondence to a neighbourhood of $\overline\Omega$. Then $f$ extends holomorphically to a (possibly smaller) neighbourhood of $\overline\Omega$.
\end{theorem}

The original proof of this result \cite{DiPi2} used the fact that $f$ extends smoothly to pseudoconvex points of $\partial\Omega$. 
This result in turn used subelliptic estimates for the $\overline\partial$-Neumann operator. Later Pinchuk-Shafikov \cite{PiSha} 
gave a self-contained geometric proof without using the $\overline\partial$-methods. Further, in~\cite{DiPi4} Diederich and 
Pinchuk showed that for holomorphic extension of the map $f$ it is enough to assume that its graph extends as an analytic set of dimension $n$ (i.e., the projection $\pi$ from this set is not assumed to be proper).

\section{Analytic discs}

In this section we consider a special case of proper holomorphic mappings from the unit disc to domains in $\C^n$. Since the unit 
disc does not have biholomorphic invariants, analytic discs are more flexible than holomorphic mappings between domains in 
$\C^n$ for $n > 1$. This flexibility makes them very useful in geometric complex analysis and its applications. These applications 
are often based on the existence of analytic discs with boundaries in prescribed CR manifolds. We discuss here some important results of this type.

\subsection{Gromov's theorem}

Consider the standard symplectic form on $\C^n$
$$
\omega = \sum_{j=1}^n dx_j \wedge dy_j.
$$
A real submanifold $E$ of dimension $n$ in $\C^n$ is called {\it Lagrangian} if $\omega\vert L = 0$. It is easy to see that every Lagrangian manifold is totally real, but the class of totally real manifolds is larger. The following result is due to Gromov \cite{Gr}.

\begin{theorem}\label{continuity10}
Let  $E$ be  a smooth compact Lagrangian submanifold  in $\cx^n$. Then there exists a nonconstant analytic  disc smooth on $\overline{\D}$ with the boundary attached to~$L$.
\end{theorem}

This theorem has deep applications in symplectic geometry (see, for example, \cite{Ar}). Note that one can view it as a 
(partial) generalization of the Riemann mapping theorem. Indeed, when $n=1$ every real curve is Lagrangian.

From the analytic point of view the problem of constructing an analytic disc with the boundary glued to $E$ can be viewed as a 
Riemann-Hilbert type boundary value problem with nonlinear boundary data (given by $E$). We sketch the main steps of 
Gromov's approach following the work of Alexander \cite{Al}, who gave a simplified version of Gromov's 
approach in the case of $\C^n$. We note that the original methods of Gromov lead to considerably more general results.

 {\it Step 1. Manifolds of discs and elliptic estimates.} Fix a point $p \in E$   and  fix also a  noninteger $r > 1$.
  Consider the set of pairs 
\begin{equation}\label{mnfd-U}
{\mathcal F} = \left\{f\in C^{r+1}(\D, \cx^n): 
f(\partial \D) \subset E,\ f(1) = p \right\}.
\end{equation}

Denote by  $F$  an open subset of ${\mathcal F}$ which consists of $f$ homotopic to a constant map
 $f^0 \equiv p$ in ${\mathcal F}$. It is well-known that $F$ is a complex Banach manifold.
Denote by $G$ the complex Banach space of all $C^r$ maps $g:\D \to \cx^n$. Set
$H = \{ (f,g) \in F \times G: \partial f /\partial \overline\zeta = g \}$.
Then $H$ is a connected submanifold of $F \times G$.  

For $0<t<1$, let $\D_t := t \D$, and $\D^+_t := t\D \cap \{\Im \zeta >0\}$.

\begin{lemma}\label{AlEst2}
 Let  $f_k:(\D_t^+,\partial \D_t^+ \cap \rl) \to (\cx^n,E)$ 
be maps of class $C^{r+1}$ that converge uniformly to $f:(\D_t^+,\partial \D_t^+ \cap \rl) \to (\cx^n,E)$. Suppose 
that the sequence $g_k = \partial f_k/\partial\overline\zeta$ converges in $C^r(\D^+_t)$ to $g \in C^r(\D^+_t)$. Then 
for every $\tau < t$ one has $f \in C^{r+1}(\D^+_t)$ and $\{ f_k \}$ converges to $f$ in $\D^+_\tau$ in the  $C^{r+1}$ norm.
\end{lemma}

Denote by $T_\D f = (2\pi i)^{-1}f * (1/ \zeta)$ the Cauchy-Green integral on $\D$. Recall the classical  regularity property of the Cauchy-Green integral: for every 
noninteger $s > 0$ the linear map $T_\D:C^s(\D) \longrightarrow C^{s+1}(\D)$ is bounded. The proof of 
Lemma~\ref{AlEst2} given in \cite{Al} is based on the standard elliptic ``bootstrapping" argument employing the above
regularity of the Cauchy-Green operator and elementary  estimates  of the harmonic measures. Notice that this proof is purely 
local, i.e., all estimates and the convergence are established in a neighbourhood of a given boundary point of a disc. The global statement is the following 

\begin{lemma}\label{AlEst1}
 Suppose that  a sequence $\{ f_k\}$ in ${\mathcal F}$ 
converges to a continuous mapping $f:(\D,\partial\D) \to (\cx^n,E)$ uniformly on 
$\overline \D$, and $g_k:= \partial f_k/\partial \overline\zeta$ converges in $C^r(\D)$ to $g \in C^r(\D)$. 
Then $f \in  C^{r+1}(\D)$ and  $\{f_k\}$ converges to $f$ 
in ${\mathcal F}$ after possibly passing to a subsequence.
\end{lemma}

Considering a finite covering of $\partial\D$ by such neighbourhoods we obtain $C^{r+1}$ convergence in a
neighbourhood of $\partial\D$. The convergence in the interior of $\D$ follows, since $f_k = T_\D g_k + h_k$, and the 
bounded sequence $\{h_k\}$ of holomorphic functions is a normal family. 

Notice that the above boundary regularity and convergence results for analytic discs are quite similar to the tools used 
in the proof of Fefferman's mapping theorem.

\medskip

{\it Step 2: Renormalization and scaling.} The canonical projection $\pi: H \to G$ given by $\pi(f,g) = g$ is a map of class $C^1$ 
between two Banach manifolds. It is known \cite{Al, Gr} that $\pi$ is  a Fredholm map of index $0$ and the 
constant map $f^0$  is a regular point for $\pi$.

The crucial property of $\pi$ is proved in \cite{Al}: {\it the map $\pi$ is not surjective.} Now, arguing by absurd, suppose
that a nonconstant analytic disc of class $C^{r+1}(\D)$ attached to $E$ does not exist, then 
$\pi^{-1}(0) = \{ f^0\}$. It follows that $0 \in G$ is a regular value of $\pi$. If $\pi$ is proper, then Gromov's argument 
based on Sard-Smale's theorem implies surjectivity of $\pi$ (see \cite{Al}) -- a contradiction. Thus, it remains to show 
that $\pi:H \to G$ is proper.

\smallskip

 Arguing by contradiction, suppose that $\pi$ is not proper. Then there exists a sequence $\{(f_k,g_k)\}\subset H$ 
 such that $g_k\to g$ in $G$ but $f_k$ diverge in $F$. For every $k$ consider the function $q_k$ defined by   
 $q_k(\zeta)= T_\D g_k(\zeta)$  for $\zeta \in \overline\D$ and $q_k(\zeta) = 0$  on $\zeta \in \cx\setminus \overline\D$.
 Then $q_k \to q = T_\D g$ in $C^{r+1}(\overline\D , \cx^n)$ and $f_k = q_k + h_k$, where 
$h_k\in C^{r+1}(\D , \cx^n)$ and $h_k$ is holomorphic on $\D$.  We have  $f_k(\partial \D)\subset E$ 
and $q_k$ are uniformly bounded since $g_k$ are; we conclude that $h_k|_{\partial \D}$ are uniformly  bounded. 
By the maximum principle the functions $h_k$ are uniformly bounded on $\overline\D$. Hence, $f_k$ are uniformly 
bounded.

Set $M_k = \sup_{\D}|h'_k(\lambda)|$. Since $h_k\in C^r(\D , \cx^n)$ and $r > 1$, 
the constants $M_k$ are finite for every~$k$. If $\{M_k\}$ contains a bounded subsequence, then a subsequence of 
$\{h_k\}$ converges uniformly on~$\D$. Then a subsequence of $\{f_k\}$ converges uniformly, and by Lemma~\ref{AlEst1}
it converges in $F$ -- a contradiction. Thus, we may  suppose that $M_k \to \infty$.  The key idea of 
\cite{Al} is to apply a renormalization argument which is essentially a version of the scaling argument.

There exists $\lambda_k\in \partial \D$ with $M_k= |h'(\lambda_k)|$ and, taking a subsequence if necessary, suppose that $\lambda_k\to \lambda^*$.
Set $z_k = (1-\frac{1}{M_k})\lambda_k\in \D$ and consider the renormalization sequence of the Mobius maps  
$$
\phi_k(\lambda) = \frac{\lambda + z_k}{1+\bar z_k\lambda} .
$$ 
Set $\tilde f_k = f_k\circ \phi_k$, 
$\tilde q_k = q_k\circ \phi_k$ and $\tilde h_k = h_k\circ \phi_k$. It is proved in \cite{Al} that after extracting 
a subsequence, the sequences $(\tilde q_k)$ and $(\tilde h_k)$ converge uniformly on compacts in 
$\overline\D\setminus \{-\lambda^*\}$ respectively to a constant map $c$  and a holomorphic map $\tilde h$.

 Notice that since $q_k$ converge in $C^{r+1} (\overline \D) $, the sequence $\tilde q_k$ converges  on compacts in 
$\overline\D\setminus \{-\lambda^*\}$ in this norm. Since Lemma~\ref{AlEst1} is local, it applies and gives the 
convergence of $(\tilde f_k)$ to $\tilde f$ also in the $C^{r+1}$-norm on compacts in 
$\overline\D\setminus \{-\lambda^*\}$. 
Then again the argument of \cite{Al} shows that $\vert\tilde h_k'(\lambda_k)\vert$ converges to  
$1/2 = \vert \tilde h'(\lambda^*) \vert$. Hence, $\tilde f$ is a nonconstant disc  of class 
$C^{r+1}$.  By the boundary regularity theorem for analytic discs, we conclude that $f$ is of class $C^\infty$
 on $\overline\D \setminus \{ 1 \}$.  

Furthermore, since $E$ is a Lagrangian manifold, it is easy to see that the disc $\tilde f$ has  
bounded area. This is due to the fact that for an analytic disc the area (induced by the Euclidean structure) 
coincides with the symplectic area $[\D] (f^*\omega)$ (where $[\D]$ is the current of integration over $\D$). 
Then \cite[Thm~2]{Al} implies that $f: \D\setminus f^{-1}(E) \to \cx^n \setminus E$ is a proper 
map. But then $\tilde f$ extends smoothly to a neighbourhood of the point $1$ (see Proposition \ref{CCS}) and 
so is smooth on $\overline\D$. This contradicts our 
assumption of nonexistence of nonconstant analytic discs attached to $E$, and  the theorem is proved.

\medskip

The assumption that $E$ is Lagrangian is crucial in Gromov's theorem. Alexander \cite{Al4} constructed a totally real 
torus $T^2$ in $\C^2$ which does not contain the boundary of an analytic disc. However, in this example one can attach 
to $T^2$ the boundary of some Riemann surface (an annulus). 
This phenomenon was recently studied by Duval-Gayet \cite{DuGa} for certain classes of totally real tori in $\C^2$. 
Their approach uses the results of Bedford-Gaveau \cite{BeGa}, Bedford-Klingenberg \cite{BeKl} and Kruzhilin~\cite{Kru2} 
on filling topological 2-spheres, contained in compact strictly pseudoconvex hypersurfaces in~$\C^2$, with Levi-flat 
hypersurfaces. The filling is provided by a 1-parameter family of analytic discs attached to the sphere.
These results have many other applications, in particular, in symplectic topology.

\subsection{Discs in pseudoconvex domains} 
Another approach to the extension of the Riemann mapping theorem concerns construction of proper holomorphic discs in 
domains in $\C^n$. We begin with the following result of Forstneri\v c -Globevnik \cite{ForGl}.

\begin{theorem}\label{FG}
Let $\Omega$ be a smoothly bounded strictly pseudoconvex domain in $\C^n$. Then for every point $p \in \Omega$ there exists an
analytic disc $f:\D \to \Omega$, smooth on $\overline{\D}$, and such that $p = f(0)$ and $f(\partial\D) \subset \partial\Omega$.
\end{theorem}

In fact, even stronger results in this direction have been obtained. We refer the reader to \cite{DrDrFor} for a detailed account.
The proof of  Theorem \ref{FG} can be described as follows. Consider a global defining function $\rho$ of $\Omega$. The idea is to construct an analytic disc attached to a suitable noncritical sub-level set of $\rho$. When such a sub-level is a small deformation of a ball around $p$, this can be achieved by the implicit function theorem. The main part of the proof consists of two major steps. First, an approximate solution of the Riemann-Hilbert type boundary value problem allows one to construct a homotopy on the space of analytic discs attached to the noncritical level sets of $\rho$. The second step is a careful analysis of the Morse geometry of a critical level set of $\rho$ which allows one to push an analytic disc on the post-critical level set. Combining these two tools  we may begin with a 
small disc attached to  some noncritical level and then deform it through other levels to a global disc attached to the boundary.

For strictly convex domains even stronger results can be obtained. This theory was developed by Lempert \cite{Lem1}. 
Let $\Omega$ be a  strictly convex domain in $\C^n$ (this means the real Hessian of the boundary is positive definite; in particular, 
$\Omega$ is strictly pseudoconvex).
Fix a point $p \in \Omega$. Then for every tangent vector $V$ at $p$ there exists a unique  analytic disc $f$ 
centred at  $p$ in the direction of $V$, which is extremal for the Kobayashi-Royden metric of $\Omega$. The  condition of extremality means that the infimum in the definition of the metric is achieved on this disc. It turns out that $f$ is smooth up to the boundary and 
its boundary is attached to $\partial\Omega$. Moreover, $f$ admits a holomorphic lift which is attached to the projectivization of the holomorphic tangent bundle of $\partial\Omega$. 
Since $\PP H(\partial\Omega)$ is a totally real manifold, $f$ satisfies a Riemann-Hilbert type boundary value problem. When $\Omega$ is a small deformation of the unit ball, this problem can be easily solved by the Implicit Function Theorem. The general case requires more advanced tools provided by the continuity method. It consists of two major steps: the Implicit Function Theorem for the linearized Riemann-Hilbert boundary value problem and a priori estimates (here the assumption of strict convexity is used).

In the case of the unit ball Lempert's discs through the origin are just linear and are given by the intersection of complex lines with the ball. It turned out that in the general case the geometry of extremal discs through any point $p$ is similar: they form a singular foliation of $\Omega$ with a unique singularity at $p$. This allows one to construct the ``Riemann mapping" from $\Omega$ to $\B^n$ which is holomorphic along every extremal disc through a fixed point $p$ and preserves the contact structure of the boundary. 

Lempert's theory has many applications. For example,  it provides an independent proof of Fefferman's mapping theorem. Furthermore, extremal discs form a very useful family of biholomorphic invariants which leads to a solution of the biholomorphic equivalence problem \cite{Lem2}. The logarithm of the Euclidean norm of the above ``Riemann mapping" gives a solution of the complex Monge-Amp\`ere equation with a logarithmic pole at $p$; it can be also viewed as a higher dimensional analog of the Green function. A similar approach was used by Donaldson \cite{Do} for the construction of regular solutions of the Dirichlet problem for a certain class of complex Monge-Amp\`ere equations.

\section{Positive codimension}

In this section we consider the properties of holomorphic mappings $f:\Omega \to \Omega'$, where $\Omega \subset \C^n$ and $\Omega' \subset \C^N$ with $1 < n < N$ (the case of {\it  positive codimension}). These maps do not have flexibility of analytic 
discs since the boundary of the source domain has  intrinsic geometry. Nevertheless, the case of positive codimension is considerably more flexible than the equidimensional one. This is illustrated by the following result  due to Forstneri\v c \cite{For1} and Low~\cite{Low}.

\begin{theorem}
Let $\Omega$ be a bounded strictly pseudoconvex domain with $C^2$ boundary in $\C^n$. There is an integer $N_1$ such that for every $N \ge N_1$ there exists a proper holomorphic mapping $f:\Omega \to \B^N$. Some of these embeddings extend continuously to $\overline\Omega$ but there exist also embeddings that are not continuous on $\partial\Omega$.
\end{theorem} 

In particular, a direct analog of  Fefferman's mapping theorem is not true in the case of positive codimension.

Note that the tools of the Moser or the  Cartan-Chern theory do not seem to be appropriate in this case. This is one of the 
reasons why the study of the rigidity phenomenon in positive codimension is a difficult problem. One of the main tools here is the geometric Reflection Principle based on the geometry of Segre varieties, which seems to admit some generalization to the case of positive codimension.

Forstneri\v c~\cite{For7} showed that most generic real analytic CR manifolds of positive CR dimension are not locally holomorphically embeddable to the germ of any generic real algebraic CR manifold of the same real codimension. 
One of the principal facts is that  an analog of  the Poincar\'e - Alexander phenomenon  holds for CR mappings between 
real spheres of positive codimension if the initial regularity of a CR mapping is sufficiently high. Forstneri\v c \cite{For2} proved that such a CR mapping extends to a rational mapping with un upper bound on the degree (depending on the codimension). Similar results are obtained for holomorphic mappings  between real algebraic CR manifolds; see for example \cite{BaEbRo3,CoMeSu,Me1,Z}. 
However, the extension of the Poincar\'e - Alexander phenomenon to the real analytic category meets difficulties. The following unpublished result is due to S.~Pinchuk (Thesis, Chelyabinsk, 1979):

\begin{theorem}
\label{hcod1}
Let $\Gamma$ be a (connected) real analytic  strictly pseudoconvex hypersurface  in $\C^n$ ($n > 1$). Assume that $f$ is a  smooth CR mapping in a neighbourhood $U$ of a point $p \in \Gamma$ such that $f(U) \subset \partial\B^N$ with $n \le N$. Then $f$ extends as holomorphic mapping along any path in $\partial\Omega$.
\end{theorem}
The proof is based on the analytic Reflection Principle. Currently there is no extension of this result to the case 
when the sphere $\partial\B^n$ is replaced with a real analytic strictly pseudoconvex hypersurface. Furthermore, even in the case of  the local Schwarz type Reflection Principle many basic questions remain open.

Finally, Forstneri\v c \cite{For3} established the following 

\begin{theorem}
\label{hcod2}
Let $f:\Gamma \to \Gamma'$ be a smooth CR mapping between real analytic strictly pseudoconvex hypersurfaces in $\C^n$ and $\C^N$ respectively, $n \le N$. Then there exists an open dense subset $O \subset \Gamma$ such that $f$ extends holomorphically 
to a neighbourhood of every point of $O$.
\end{theorem}

A natural question is if the above set $O$ coincides with the whole $\Gamma$. The following result was obtained in \cite{PiSu}.

\begin{theorem}
\label{hcod3}
Let $f:\Gamma \to \Gamma'$ be a smooth CR mapping between real analytic strictly pseudoconvex hypersurfaces in $\C^n$ and $\C^N$ respectively, and 
$n \le N \le 2n$. Then $f$ extends holomorphically to a neighbourhood of every point in $\Gamma$.
\end{theorem}

The proofs of Theorems \ref{hcod2} and \ref{hcod3} are based on the geometric Reflection Principle and the study of Segre 
varieties. To the best of our knowledge it is not known whether the condition $N \le 2n$ in Theorem \ref{hcod3} can be dropped.


\end{document}